\documentclass{article}

\usepackage{amssymb}
\usepackage[leqno]{amsmath}
\usepackage{amsfonts}
\usepackage{amsopn}
\usepackage{amstext}
\usepackage{amsthm}

\providecommand{\cal}{\mathcal}
\renewcommand{\Bbb}{\mathbb}

\newenvironment{pf}{\begin{proof}}{\end{proof}}

\newcommand{\Aa}{{\Bbb{A}}}

\newcommand{\Ef}{{\cal{F}}}

\newcommand{\Ai}{{\cal{I}}}

\newcommand{\Pee}{{\cal{P}}}

\newcommand{\Tau}{{\tau}}

\newcommand{\Be}{{\Bbb{B}}}
\newcommand{\Nat}{{\Bbb{N}}}

\newcommand{\Err}{{\Bbb{R}}}

\newcommand{\lam}{{\lambda}}
\newcommand{\al}{\alpha}
\newcommand{\Gam}{\Gamma}

\newcommand{\sig}{\sigma}
\newcommand{\eps}{\varepsilon}
\renewcommand{\phi}{\varphi}
\renewcommand{\rho}{\varrho}

\newcommand{\rest}{\restriction}

\newcommand{\ntr}{n\in\omega}

\newcommand{\loe}{\leqslant}
\newcommand{\goe}{\geqslant}

\newcommand{\subs}{\subseteq}
\newcommand{\sups}{\supseteq}
\newcommand{\nnempty}{\ne\emptyset}
\newcommand{\argum}{\,\cdot\,}
\newcommand{\ovr}{\overline}

\newcommand{\cl}{\operatorname{cl}}

\newcommand{\dist}{\operatorname{dist}}

\newcommand{\id}{\operatorname{id}}

\newcommand{\liminv}{\varprojlim}
\newcommand{\h}{\widehat}

\newcommand{\Es}{{\cal{S}}}

\newcommand{\Land}{\;\&\;}

\newcommand{\Fin}{\operatorname{fin}}

\newtheorem{tw}{Theorem}[section]
\newtheorem{wn}[tw]{Corollary}
\newtheorem{lm}[tw]{Lemma}
\newtheorem{prop}[tw]{Proposition}
\newtheorem{claim}[tw]{Claim}
\theoremstyle{definition}

\newtheorem{ex}[tw]{Example}
\theoremstyle{remark}
\newtheorem{uwgi}[tw]{Remark}

\newcommand{\setof}[2]{\{#1\colon #2\}}
\newcommand{\bigsetof}[2]{\Bigl\{#1\colon #2\Bigr\}}

\newcommand{\sett}[2]{\{#1\}_{#2}}
\newcommand{\sn}[1]{\{#1\}} % singleton
\newcommand{\map}[3]{#1\colon #2 \to #3} % A function
\newcommand{\img}[2]{#1[#2]} % image of a set
\newcommand{\dpower}[2]{[#1]^{#2}}

\newcommand{\fin}[1]{[#1]^{<\omega}}

\newcommand{\Sig}{\Sigma}
\newcommand{\Ssig}{{\sig\in\Sigma}}

\newcommand{\I}{\ensuremath{\mathcal I}}
\providecommand{\nat}{\omega}

\newcommand{\im}{\operatorname{im}}
\newcommand{\bal}{\operatorname{B}}
\newcommand{\clbal}{\overline{\bal}}
\newcommand{\anorm}{\|\cdot\|}
\newcommand{\norm}[1]{\|#1\|}
\newcommand{\bnorm}[1]{\Bigl\|#1\Bigr\|}

\newcommand{\aabs}{|\cdot|}
\newcommand{\abs}[1]{|#1|}
\newcommand{\uball}[1]{\clbal_{#1}}
\newcommand{\ubal}{\uball}
\newcommand{\usphere}[1]{\operatorname{S}_{#1}}
\renewcommand{\S}{{\mathbb S}}
\newcommand{\R}{\ensuremath{\mathcal R}}
\newcommand{\N}{\ensuremath{\mathbb N}}
\newcommand{\ult}[1]{\operatorname{Ult}({#1})}

\title{Kadec norms on spaces of continuous functions}
\author{
{\sc Maxim R.~Burke} \thanks{Research supported by NSERC. The
first and second authors thanks J.~Steprans and the organizers of
the Thematic Program on Set Theory and Analysis at The Fields
Institute for Research in Mathematical Sciences for helping to
arrange the second author's postdoctoral fellowship at the
University of Prince Edward Island from January to April 2003.}
\\
{\footnotesize Department of Mathematics and Statistics} \\
{\footnotesize University of Prince Edward Island} \\
{\footnotesize Charlottetown PE, Canada C1A 4P3} \\
{\footnotesize burke@upei.ca}
\and {\sc Wies{\l}aw Kubi\'s} \thanks{Research supported by the
NSERC grant of the first author.}
\\
{\footnotesize Department of Mathematics and Statistics} \\
{\footnotesize University of Prince Edward Island} \\
{\footnotesize Charlottetown PE, Canada C1A 4P3} \\
{\footnotesize kubis@ux2.math.us.edu.pl}
\and {\sc Stevo Todor\v{c}evi\'c}
\thanks{The third author thanks the Department of
Mathematics and Statistics at the University of Prince Edward
Island for its hospitality during his visits while this paper was
being written.
\endgraf
AMS Subject Classification. Primary: 46B03, 46B26. Secondary: 46E15, 54C35.
\endgraf Key words and phrases: $\tau_p$-Kadec norm, Banach space of
continuous functions, compact space.
 }
\\
{\footnotesize Universite Paris 7} \\
{\footnotesize CNRS, UMR 7056} \\
{\footnotesize 2, Place Jussieu} \\
{\footnotesize 75251 Paris Cedex 05} \\
{\footnotesize France} \\
{\footnotesize stevo@logique.jussieu.fr}
 }

\begin{document}
\maketitle

\begin{abstract}
We study the existence of pointwise Kadec renormings for Banach
spaces of the form $C(K)$. We show in particular that such a
renorming exists when $K$ is any product of compact linearly
ordered spaces, extending the result for a single factor due to
Haydon, Jayne, Namioka and Rogers. We show that if $C(K_1)$ has a
pointwise Kadec renorming and $K_2$ belongs to the class of spaces
obtained by closing the class of compact metrizable spaces under
inverse limits of transfinite continuous sequences of retractions,
then $C(K_1\times K_2)$ has a pointwise Kadec renorming. We also
prove a version of the three-space property for such renormings.
\end{abstract}

\section{Introduction}

Let $X$ be a Banach space. Let $\Tau$ be a tvs topology on $X$
weaker than the norm topology. The norm on $X$ is called {\em
$\Tau$-Kadec} if the norm topology coincides with $\Tau$ on the
unit sphere. When $\Tau$ is the weak topology, the norm is simply
said to be Kadec. In our setting we consider mainly spaces of the
form $X=C(K)$ for some compact space $K$. We shall be interested
primarily in the question of when there is a norm on $X$
equivalent to the supremum norm which is $\Tau_p$-Kadec where
$\Tau_p$ stands for the topology of pointwise convergence,
referred to henceforth as the pointwise topology.

Raja has shown in \cite{Ra} that the existence of a $\Tau$-Kadec
renorming for $X$ is equivalent to the existence of a countable
collection $\{A_n:n\in\N\}$ of convex subsets of $X$ such that the
collection of sets of the form $U\cap A_n$, where $U\in \Tau$,
forms a network for the norm topology.  (A collection $C$ of sets
in a topological space is a {\em network} for the topology if
every open set is the union of a subcollection of $C$. In other
words, $C$ is like a base except that its members do not have to
be open.) It is not known whether the word ``convex'' can be
omitted in this characterization. The notion obtained by deleting
convexity goes by several names in the literature. Following
\cite{JNR1} (where the notion was introduced), we say that {\em
$(X,\Tau)$ has a countable cover by sets of small local
norm-diameter}, or more briefly {\em $(X,\Tau)$ is norm-SLD}, if
there is a countable collection $\{A_n:n\in\N\}$ of subsets of $X$
such that the sets $U\cap A_n$, where $n\in\N$ and $U\in \Tau$,
form a network for the norm topology. It is shown in \cite{JNR2}
that when $K$ is an infinite compact $F$-space, then $C(K)$ is not
$\sigma$-fragmentable, in particular $C(K)$ has no Kadec
renorming.

In the paper \cite{HJNR}, it is shown that for every compact
totally ordered space $K$, $C(K)$ has a $\Tau_p$-Kadec renorming.
We shall show that the conclusion remains true if $K$ is an
arbitrary product of compact linearly ordered spaces. This
improves the result in \cite[Theorem~5.21(b)]{Bu1} (due to Jayne,
Namioka and Rogers for countable products, see \cite[Remark~(1),
p.~329]{JNR3}) that for such a product $K$, $C(K)$ is norm-SLD in
the pointwise topology. It is unknown whether the existence of a
$\Tau_p$-Kadec renorming for each of $C(K_1)$ and $C(K_2)$ implies
the existence of such a renorming for $C(K_1\times K_2)$. Ribarska
has shown in \cite{Ri2} that if $C(K_1)$ has a $\Tau_p$-Kadec
renorming and $C(K_2)$ is norm-SLD in the pointwise topology, then
$C(K_1\times K_2)$ is norm-SLD in the pointwise topology. We
establish that if $C(K_1)$ has a $\Tau_p$-Kadec renorming and
$K_2$ belongs to the class of spaces obtained by closing the class
of compact metrizable spaces under inverse limits of transfinite
continuous sequences of retractions, then $C(K_1\times K_2)$ has a
$\Tau_p$-Kadec renorming.

In \cite{LZ}, the authors establish, under certain conditions, the
three-space property for a sequential version of the Kadec
property. (A property of Banach spaces is a three-space property
if $X$ has the property whenever $Y$ and $X/Y$ do, where $Y$ is a
subspace of $X$.) A Banach space is said to have the {\em
Kadec-Klee property} if every weakly convergent sequence on the
unit sphere is strongly convergent. (The terminology is not used
consistently in the literature. In particular, in \cite{DGZ} a
norm which has the Kadec-Klee property is what we have called a
Kadec norm.) A norm is {\em locally uniformly rotund} (LUR) if
whenever $x_n$, $n\in\N$, and $x$ are on the unit sphere and
$\lim\|x_n+x\|=2$ we have $\lim x_n=x$. As pointed out in
\cite{Al}, if the norm in a Banach space $X$ is LUR and $\Tau$ is
a tvs topology on $X$ such that the unit ball is $\Tau$-closed
(for example the weak topology), then the norm is necessarily
$\Tau$-Kadec. In \cite{LZ}, it is shown that if $X$ is a Banach
space, $Y$ is a subspace of $X$, $Y$ has the Kadec-Klee property
and $X/Y$ has an LUR renorming, then $X$ has the Kadec-Klee
property. We show, solving a problem raised in \cite{LZ}, that the
Kadec-Klee property can be replaced by the Kadec property in their
result. It is not known whether the existence of a Kadec renoming
is a three-space property. Ribarska has shown in \cite{Ri1} that
being norm-SLD in the weak topology is a three-space property. Her
proof also shows that for spaces $L\subs K$, if $C(L)$ and
$C_0(K\setminus L)$ are norm-SLD in the pointwise topology, then
so is $C(K)$.

We write lsc, usc for lower semi-continuous, upper
semi-continuous, respectively. Given a map $\map fXY$, a {\em
level set} of $f$ is any set of the form $\setof{x\in
X}{f(x)=y_0}$, where $y_0\in Y$ is fixed. Given a normed space
$(X,\anorm)$ we denote by $\ubal X$ and $\usphere X$ the closed
unit ball and the unit sphere of $X$ respectively. A closed (resp.
open) ball centered at $x$ and with radius $r>0$ is denoted by
$\clbal(x,r)$ (resp. $\bal(x,r)$). Similarly, for a set $A\subs
X$, $\bal(A,r)$ denotes $\setof{x\in
X}{\dist(x,A)<r}=A+\bal(0,r)$.

\section{Preliminaries}

We begin with a standard fact.

\begin{prop}\label{continuous image}
Let $K$ and $L$ be compact spaces, and let $\phi\colon K\to L$ be
a continuous surjection. Then the map $T\colon C(L)\to C(K)$
defined by $T(f)=f\phi$ is a linear isometry and a
$\tau_p$-homeomorphism onto its range. In particular, if $C(K)$
has an equivalent $\tau_p$-Kadec norm, then so does $C(L)$.
\end{prop}

\begin{pf}
$T$ is clearly linear. We have
$\|T(f)\|_\infty=\|f\phi\|_\infty=\|f\|_\infty$ because $\phi$ is
onto, so $T$ is an isometry. The fact that $T$ is a
$\tau_p$-homeomorphism onto its range follows from the fact that
$\phi$ is onto and from the equality $T(f)(x)=f(\phi x)$ for $x\in
K$ .
\end{pf}

The following Proposition is given as \cite[Proposition~1]{Al} for
the case where $\Tau$ is generated by a total subspace of $X^*$.
As pointed out in \cite[Proposition 4]{Ra}, the proof works for
any linear topology.

\begin{prop}\label{alexandrov}
Let $X$ be a Banach space whose norm is $\Tau$-Kadec. Then the
norm is $\Tau$-lsc, i.e., the unit ball is $\Tau$-closed.
\hspace{\stretch{1}}$\square$
\end{prop}

\begin{prop}\label{annulus} {\rm (Cf.~\cite[Lemma~1]{Ra}.)}
Let $X$ be a Banach space, $x_0\in S_X$, $\Tau$ a weaker linear
topology on $X$ with respect to which the norm is $\Tau$-Kadec at
$x_0$ {\rm (}i.e., the norm and $\Tau$ neighborhoods of $x_0$ are
the same{\rm )}. Then for any $r>0$, there exists $\delta>0$ and a
neighborhood $U\in\tau$ of $x_0$ such that
$$
U\cap\bal(0,1+\delta)\subs \bal(x_0,r).
$$
\end{prop}

\begin{pf}
Find a neighborhood $W\in\tau$ of $x_0$ such that $W\cap \usphere
X\subs \bal(x_0,r/2)$. By the $\tau$-continuity of the addition,
there are $V,V'\in\tau$ such that $x_0\in V$, $0\in V'$ and
$V+V'\subs W$. Fix $\delta>0$ such that $\delta\loe r/2$ and
$\bal(0,\delta)\subs V'$. Then $V\cap (\usphere
X+\bal(0,\delta))\subs \bal(x_0,r)$. Indeed, if $y\in V$ and
$\norm{y-z}<\delta$ for some $z\in \usphere X$ then $z\in
(V+V')\cap\usphere X\subs\bal(x_0,r/2)$ so
$\norm{y-x_0}\loe\norm{y-z}+\norm{z-x_0}<r/2+\delta\loe r$. As
closed balls are $\tau$-closed (Proposition~\ref{alexandrov}), we
may assume that $V\cap\clbal(0,1-\delta)=\emptyset$. Then $V\cap
\bal(0,1+\delta)\subs\bal(x_0,r)$.
\end{pf}

We shall need the simple facts about lower semi-continuous maps
given by the next three propositions and their corollaries.

\begin{prop}\label{snieg0} Let $X$ be a topological space and let
$\map{f,g}X{\Err}$ be functions whose sum is identically equal to
a constant value $k\in\Err$. For any $x\in X$, if $f$ is lsc at
$x$, then $g$ is usc at $x$.
\end{prop}

\begin{pf}
Fix $\eps>0$ and find a neighborhood $V$ of $x$ such that
$f(x')>f(x)-\eps$ for $x'\in V$. Thus
$g(x')=k-f(x')<k-f(x)+\eps=g(x)+\eps$ whenever $x'\in V$.
\end{pf}

\begin{wn}\label{snieg1}
Let $X$ be a topological space.
\begin{enumerate}
\item[{\rm(a)}]
If $\map{f,g}X{\Err}$ are lsc, then the restrictions of $f$ and
$g$ to any level set for $f+g$ are continuous.
\item[{\rm(b)}]
If $\map{f_n}X{\Err}$, $n\in\N$, are nonnegative lsc functions
such that $\sum_{n\in\N}f_n$ converges pointwise, then the
restriction of each $f_n$ to a level set for $\sum_{n\in\N}f_n$ is
continuous.
\end{enumerate}
\end{wn}

\begin{pf} (a) Applying Lemma \ref{snieg0} to the restrictions of $f$ and
$g$ to a level set $S=\{x\in X:f(x)+g(x)=k\}$ shows that because
these functions are lsc at every point, they are also usc at every
point.

(b) Apply part (a) to $f=f_n$ and $g=\sum_{m\not=n}f_m$.
\end{pf}

It will be useful to have a slightly stronger version of
Corollary~\ref{snieg1}(b).

\begin{prop}\label{snieg2}
Let $X$ be a topological space, $\sett{f_n}{\ntr}$ a sequence of
nonnegative lsc real-valued functions on $X$ such that
$\theta(x)=\sum_{\ntr}f_n(x)$ is finite for every $x\in X$. Assume
$\sett{x_\sig}{\Ssig}$ is a net in $X$ converging to $x\in X$  and
$\lim_{\Ssig}\theta(x_\sig)=\theta(x)$ for every $\Ssig$. Then
$\lim_{\Ssig}f_k(x_\sig)=f_k(x)$ for every $k\in\nat$.
\end{prop}

\begin{pf}
Fix $k\in\nat$ and let $g=\sum_{n\ne k}f_n$. Observe that $g$ is
lsc as the supremum of a set of lsc functions. Fix $\eps>0$. There
exists $\sig_0$ such that $\theta(x_\sig)-\theta(x)\loe\eps/2$,
$f_k(x_\sig)>f_k(x)-\eps$ and $g(x_\sig)>g(x)-\eps/2$ for
$\sig\goe\sig_0$. Fix $\sig\goe \sig_0$ and suppose
$f_k(x_\sig)\not<f_k(x)+\eps$. Then
$$
\theta(x_\sig)=f_k(x_\sig)+g(x_\sig) >
f_k(x)+\eps+g(x)-\eps/2=\theta(x)+\eps/2,
$$
so $\theta(x_\sig)-\theta(x)>\eps/2$, a contradiction.
\end{pf}

\begin{prop}\label{lsc1}
Let $X$ be a topological space, $n\in\N$, $f_i\colon X\to\Err$
for $1\leq i\leq n$. Let $x\in X$. Suppose $\sum f_i\leq 0$, $\sum
f_i(x)=0$, and each $f_i$ is lsc at $x$. Then each $f_i$ is
continuous at $x$.
\end{prop}

\begin{pf}
Fix $i$ and $\eps>0$. For $y$ in some neighborhood of $x$ we
have
\[
{\textstyle f_i(x)-\eps<f_i(y)\leq -\sum_{j\not= i}
f_j(y)<-\left(\sum_{j\not= i}
(f_j(x)-\eps/(n-1))\right)=f_i(x)+\eps. }
\]
\end{pf}

\begin{wn}\label{lsc2}
Let $X$ be a topological space, $n\in\N$, $f_i\colon X\to\Err$ for
$1\leq i\leq n$, $h\colon X\to\Err$. Let $x\in X$. Suppose $\sum
f_i\leq h$, $\sum f_i(x)=h(x)$, each $f_i$ is lsc at $x$ and $h$
is usc at $x$. Then $h$ and each $f_i$ is continuous at $x$.
\end{wn}

\begin{pf}
$f_0+\dots+f_{n-1}-h\leq 0$,
$f_0(x)+\dots+f_{n-1}(x)-h(x)=0$ and $-h$ is lsc at $x$.
\end{pf}

An {\em inverse sequence} is a family of mappings
$p^\beta_\al\colon X_\beta\to X_\al$, $\al<\beta<\kappa$, where
$\kappa$ is a limit ordinal, such that $\al<\beta<\gamma\implies
p^\beta_\al p^\gamma_\beta = p^\gamma_\beta$. Usually, the maps
$p^\beta_\al$ are surjections. We refer the reader to
\cite[Section~2.5]{En} for the basic properties of inverse
systems. We recall here some of the relevant terminology.

We write $\S=\setof{X_\al;p^\beta_\al}{\al<\beta<\kappa}$ and we
call $p^\beta_\al$'s the {\em bonding mappings} of $\S$. The {\em
inverse limit} of $\S$, denoted by $\liminv\S$ is defined to be
the subspace of the product $\prod_{\al<\kappa}X_\al$ consisting
of all $x$ such that $p^\beta_\al(x(\beta))=x(\al)$ for every
$\al<\beta<\kappa$. If each $X_\al$ is compact then $\liminv
\S\nnempty$. If moreover each $p^\beta_\al$ is a surjection then
the projection $\map {p_\al}{\liminv\S}{X_\al}$ is also a
surjection. From a category-theoretic perspective, the inverse
limit of $\S$ is a space $X$ together with a family of continuous
maps (called {\em projections}) $\setof{p_\al}{\al<\kappa}$ which
has the property that for every space $Y$ and a family of
continuous maps $\setof{f_\al}{\al<\kappa}$ such that $p^\beta_\al
f_\beta = f_\al$ holds for every $\al<\beta<\kappa$, there exists
a unique continuous map $\map hYX$ such that $p_\al h = f_\al$ for
every $\al<\kappa$. The limit is uniquely determined in the sense
that if $X'$ with projections $p'_\al$, $\al<\kappa$, is another,
then the unique continuous map $h\colon X'\to X$ such that $p_\al
h=p'_\al$ for all $\al<\kappa$ is a homeomorphism. The definition
of $\liminv\S$ given above is one of the possibilities. We will
use the property that $\liminv\setof{X_\al;
p^\beta_\al}{\al<\beta<\kappa}$ is isomorphic to $\liminv
\setof{X_\al; p^\beta_\al}{\al<\beta,\; \al,\beta\in C}$ for every
cofinal set $C\subs \kappa$.

An inverse sequence
$\S=\setof{X_\al;p^\beta_\al}{\al<\beta<\kappa}$ is {\em
continuous} if for every limit ordinal $\delta<\kappa$ the space
$X_\delta$ together with $\setof{p^\delta_\al}{\al<\delta}$ is
homeomorphic to $\liminv\setof{X_\al;
p^\beta_\al}{\al<\beta<\delta}$.

A {\em retraction} is a continuous map $\map fXY$ which has a
right inverse, i.e. a continuous map $\map jYX$ with $f j=\id_Y$.
Note that $j$ is an embedding and $f$ restricted to $\img jY$ is a
homeomorphism.

Finally, we point out that many of our results about Banach spaces
equipped with a weaker linear topology $\Tau$ with respect to
which the norm is lsc have conclusions which assert the existence
of an equivalent norm with a certain property.  In all such
results, the assumption that the norm is $\Tau$-lsc can be
weakened to the assumption that the $\Tau$-closure of the unit
ball is bounded, since the Minkowski functional of this closure
provides an equivalent $\Tau$-lsc norm.

\section{Finite products of linearly ordered spaces}

In this section we show that $C(L_0\times\dots\times L_{n-1})$ has
a $\Tau_p$-Kadec renorming, whenever $L_0,\dots,L_{n-1}$ are
compact linearly ordered spaces. In Theorem~\ref{products of lo},
this result will be extended to arbitrary products.

\begin{lm}\label{mewa} If $X$ is a compact linearly ordered space,
$(Y,d)$ is a metric space, $\map fXY$ is continuous, and for each
$m\in\nat$ we set
$$
v_m(f)=\sup\bigsetof{ \sum_{i<m}d(f(a_i),f(a_{i+1}))}{a_0\loe
a_1\loe \dots\loe a_m},
$$
where $0$ and $1$ denote the first and last elements of $X$, then
$$\lim_{m\to\infty}\ v_{m+1}(f)-v_m(f)=0.$$
\end{lm}

\begin{pf}
Fix $\eps>0$. Let $\I$ be a finite cover of $X$ by open intervals
$I$ such that $f[I]$ has diameter $<\eps$. Fix any $m>|\I|$.
Choose $a_0\loe a_1\loe \dots\loe a_m\loe a_{m+1}$ so that
$v_{m+1}(f)=\sum_{i<m+1}d(f(a_i),f(a_{i+1}))$. For some $I\in \I$
and $i_0<m+1$, we have $a_{i_0},a_{i_0+1}\in I$. Suppose $i_0<m$.
Then
\begin{eqnarray}
d(f(a_{i_0}),f(a_{i_0+1}))+d(f(a_{i_0+1}),f(a_{i_0+2}))  &\loe&
d(f(a_{i_0}),f(a_{i_0+2})) +
2d(f(a_{i_0}),f(a_{i_0+1}))\nonumber\\
&<& d(f(a_{i_0}),f(a_{i_0+2})) + 2\eps\nonumber
\end{eqnarray}
and we get
\begin{eqnarray}
v_m(f) & \geq & d(f(a_0),f(a_1)) + \dots +
d(f(a_{i_0-1}),f(a_{i_0}))
 + d(f(a_{i_0}),f(a_{i_0+2})) \nonumber\\
& & + d(f(a_{i_0+2}),f(a_{i_0+3})) + \dots +
d(f(a_m),f(a_{m+1}))\nonumber\\
& > & \sum_{i<m+1}d(f(a_i),f(a_{i+1}))-2\eps\nonumber\\
& = & v_{m+1}(f)-2\eps,\nonumber
\end{eqnarray}
which gives $0\loe v_{m+1}(f)-v_m(f)<2\eps$. If $i_0=m$, replace
the triple $(f(a_{i_0}),f(a_{i_0+1}),f(a_{i_0+2}))$ by the triple
$(f(a_{i_0-1}),f(a_{i_0}),f(a_{i_0+1}))$ in the argument above.
\end{pf}

Let $L$ be a linearly ordered space. We say that points $x,y\in L$
are {\em adjacent} if $x\ne y$ and no point is strictly between
$x,y$.

\begin{tw}\label{wiatr}
Assume $L_i$, $i<n$ are compact linearly ordered spaces and
$D_i\subs L_i$ is dense in $L_i$ and contains all pairs of
adjacent points for each $i<n$. Then $C(\prod_{i<n}L_i)$ has an
equivalent $\Tau_p(D)$-Kadec norm, where $D=\prod_{i<n}D_i$.
\end{tw}

(See Theorem~\ref{products of lo} for the case of arbitrary
products.)

\begin{pf}
For $f\in C(\prod_{i<n}L_i)$, we will need to consider expressions
of the form
\begin{equation}
f(x_0,x_1,\dots,x_{k-1},\,a,\,x_{k+1},\dots,x_{n-1}). \tag{1}
\end{equation}
For notational convenience, we sometimes permute the arguments so
that $a$ comes first. Letting $h_k\colon L_k\times
\prod_{\ell<n,\,\ell\ne k}L_{\ell}\to \prod_{\ell<n}L_{\ell}$ be
given by
$$h_k(a,\,x_0,x_1,\dots,x_{k-1},x_{k+1},\dots,x_{n-1}) =
(x_0,x_1,\dots,x_{k-1},\,a,\,x_{k+1},\dots,x_{n-1}),$$
we can then write
$$f(h_k(a,\,x_0,x_1,\dots,x_{k-1},x_{k+1},\dots,x_{n-1}))$$
instead of (1).

For each $k<n$ and $m\in\nat$, define $v^k_m(f)$ on
$C(\prod_{i<n}L_i)$ by letting
$$
v^k_m(f)=\sup\bigsetof{
\sum_{i<m}\norm{f(h_k(a^k_i,\,\cdot\,))-f(h_k(a^k_{i+1},\,\cdot\,))}_\infty}
{
%\{a^k_0,a^k_1,\dots,a^k_m\}\subs D_k,\, a^k_0\loe
a^k_1\loe\dots\loe a^k_m}.
$$
The function $v^k_m$ is a $\Tau_p(D)$-lsc seminorm and
$$\lim_{m\to\infty}\ v^k_{m+1}(f)-v^k_m(f)=0,$$
by Lemma \ref{mewa}.

Define $\aabs$ on $C(\prod_{i<n}L_i)$ as follows.
$$
|f|=\norm f_\infty + \sum_{k<n}\sum_{m\in\nat}\frac{1}{m\cdot 2^m}v^k_m(f).
$$
It is readily seen that $\aabs$ is a norm on $C(\prod_{i<n}L_i)$
and is equivalent to the sup norm.  We now verify that it is a
$\Tau_p(D)$-Kadec norm. Since the terms in the definition of $\abs
f$ are all $\Tau_p(D)$-lsc functions of $f$,
Corollary~\ref{snieg1}(b) implies that they are all
$\Tau_p(D)$-continuous functions of $f$ when restricted to
$S:=\setof{f}{\abs f =1}$. Fix $f\in S$ and $\eps>0$.

For each $k<n$, the map $x\mapsto f(h_k(x,\,\cdot\,))$ is
continuous (with the norm topology on the range), so there is a
finite collection $\I_k$ of open intervals covering $L_k$ such
that the diameter in $C(\prod_{\ell<n,\,\ell\not=k}L_{\ell})$ of
$\setof{f(h_k(x,\,\cdot\,))}{x\in I}$ is less than $\eps$ for each
$I\in\I_k$. We may assume that $\inf I\in D_k\cup\sn0$ and $\sup
I\in D_k\cup\sn1$ for each $I\in \I_k$. Let $A_k=\setof{\inf
I}{I\in\I_k}\cup\setof{\sup I}{I\in\I_k}$. Then $A_k\subs
D_k\cup\{0,1\}$.

Let $m\in\nat$ be such that for each $k<n$,
$v^k_{m+3}(f)-v^k_m(f)<\eps$.

For each $k<n$, fix $a^k_0\loe a^k_1\loe \dots\loe a^k_m$ in $D_k$
such that
$$
v^k_m(f)<\sum_{i<m}\norm{f(h(a^k_i,\,\cdot\,))-
f(h(a^k_{i+1},\,\cdot\,))}_\infty + \delta,
$$
where $\delta=\eps/(m+4)$. Let $H_k=\setof{a^k_i}{i\loe m}\cup
(A_k\cap D_k)$.

Fix a $\Tau_p(D)$-open neighborhood $U$ of $f$ such that for $g\in
S\cap U$ we have, for all $k<n$, that $v^k_{m+i}(g)$ is strictly
within $\eps$ of $v^k_{m+i}(f)$ for $i\loe 3$. This gives
$$
|v^k_{m+i}(g)-v^k_{m+i}(f)|<\eps\quad\text{ and } \quad
|v^k_{m+i}(f)-v^k_{m}(f)|<\eps
$$
and hence
$$|v^k_{m+i}(g)-v^k_{m}(f)|< 2\eps.$$

For each $k<n$ and for each pair of elements $a<b$ of $D_k$,
choose $x=x^k_{a,b}$ and $y=y^k_{a,b}$ in $D$
such that $x(k)=a$, $y(k)=b$, $x(\ell)=y(\ell)$ for all
$\ell\not=k$ and
$$
\norm{f(h_k(a,\,\cdot\,))-f(h_k(b,\,\cdot\,))}_\infty <
|f(x)-f(y)|+\delta.
$$
Write
$$
\ovr{H}_k=H_k\cup\setof{z(k)}{z=x^{\ell}_{a,b}\ \mbox{or}\
z=y^{\ell}_{a,b}\ \mbox{for some}\ \ell<n\ \mbox{and some}\ a<b\
\mbox{in}\ H_{\ell}}.
$$
Then $\ovr{H}_k\subs D_k$.
Let $g\in U$ agree sufficiently closely with $f$ on
$H=\prod_{k<n}\ovr{H}_k$ so that $|g(h)-f(h)|<\eps$ for each $h\in
H$ and the following condition is satisfied.
\begin{enumerate}
\item[($*$)]
For each $k<n$, for each $i_0<m$, and any choice of elements of
$H_k$ of the form
$$
a^k_{i_0}=b_0\loe b_1\loe b_2\loe b_3=a^k_{i_0+1}
$$
we have, for each $j_0<3$,
\begin{eqnarray}
& &
\sum_{i\not=i_0}|g(x^k_{a^k_i,a^k_{i+1}})-g(y^k_{a^k_i,a^k_{i+1}})|
+ \sum_{j\not=j_0}|g(x^k_{b_j,b_{j+1}})-g(y^k_{b_j,b_{j+1}})|
\nonumber\\
& > &
\sum_{i\not=i_0}|f(x^k_{a^k_i,a^k_{i+1}})-f(y^k_{a^k_i,a^k_{i+1}})|
+ \sum_{j\not=j_0}|f(x^k_{b_j,b_{j+1}})-f(y^k_{b_j,b_{j+1}})| -
\eps. \nonumber
\end{eqnarray}
\end{enumerate}
Assume also that for each $k<n$ we have
\begin{equation}
\sum_{i<m}|g(x^k_{a^k_i,a^k_{i+1}})-g(y^k_{a^k_i,a^k_{i+1}})|
>
\sum_{i<m}|f(x^k_{a^k_i,a^k_{i+1}})-f(y^k_{a^k_i,a^k_{i+1}})|
-\eps. \tag{$*_1$}
\end{equation}

From ($*$) it follows that for any $x\in[b_{j_0},b_{j_0+1}]$,
writing
$$s=\norm{f(h_k(b_{j_0},\,\cdot\,))-f(h_k(x,\,\cdot\,))}_\infty +
\norm{f(h_k(x,\,\cdot\,))-f(h_k(b_{j_0+1},\,\cdot\,))}_\infty$$
and
$$t=\norm{g(h_k(b_{j_0},\,\cdot\,))-g(h_k(x,\,\cdot\,))}_\infty +
\norm{g(h_k(x,\,\cdot\,))-g(h_k(b_{j_0+1},\,\cdot\,))}_\infty$$
we have
\begin{eqnarray}
v^k_m(f) + 2\eps - t & > & v^k_{m+3}(g) - t \nonumber\\
& \goe &
\sum_{i\not=i_0}\norm{g(h_k(a^k_i,\,\cdot\,))-
   g(h_k(a^k_{i+1},\,\cdot\,))}_\infty \nonumber\\
& &+
\sum_{j\not=j_0}\norm{g(h_k(b_j,\,\cdot\,))-
   g(h_k(b_{j+1},\,\cdot\,))}_\infty
\nonumber\\
& \goe &
\sum_{i\not=i_0}|g(x^k_{a^k_i,a^k_{i+1}})-g(y^k_{a^k_i,a^k_{i+1}})|
+ \sum_{j\not=j_0}|g(x^k_{b_j,b_{j+1}})-g(y^k_{b_j,b_{j+1}})|
\nonumber\\
& > &
\sum_{i\not=i_0}|f(x^k_{a^k_i,a^k_{i+1}})-f(y^k_{a^k_i,a^k_{i+1}})|
+ \sum_{j\not=j_0}|f(x^k_{b_j,b_{j+1}})-f(y^k_{b_j,b_{j+1}})| -
\eps \nonumber \\
& > &
\sum_{i\not=i_0} (\norm{f(h_k(a^k_i,\,\cdot\,))-
   f(h_k(a^k_{i+1},\,\cdot\,))}_\infty-\delta)\nonumber\\
& &+
\left(\ \sum_{j\not=j_0}(\norm{f(h_k(b_j,\,\cdot\,))-
   f(h_k(b_{j+1},\,\cdot\,))}_\infty -\delta)
+ s\right) - s -\eps \nonumber\\
& \goe &
\sum_{i<m}\norm{f(h_k(a^k_i,\,\cdot\,))-
   f(h_k(a^k_{i+1},\,\cdot\,))}_\infty
- s - \eps - (m+3)\delta
\nonumber\\
& > & v^k_m(f)-s-2\eps\nonumber
\end{eqnarray}
and hence $t<s+4\eps$, i.e., for any $x\in[b_{j_0},b_{j_0+1}]$,
\begin{enumerate}
\item[($**$)]
$\norm{g(h_k(b_{j_0},\,\cdot\,))-g(h_k(x,\,\cdot\,))}_\infty +
\norm{g(h_k(x,\,\cdot\,))-g(h_k(b_{j_0+1},\,\cdot\,))}_\infty$

$<\norm{f(h_k(b_{j_0},\,\cdot\,))-f(h_k(x,\,\cdot\,))}_\infty +
\norm{f(h_k(x,\,\cdot\,))-f(h_k(b_{j_0+1},\,\cdot\,))}_\infty + 4\eps.$
\end{enumerate}

Consider a point $p\in \prod_{k<n}L_k$. Define
$$T=\setof{k<n}{p_k\not\in\ovr{H}_k}.$$
We will show by induction on $r=|T|$ that
$|g(p)-f(p)|<(7r+1)\eps$. This is true if $r=0$ since then $p\in
H$. For the inductive step, suppose $|T|=r+1$. Choose an open
neighborhood of $p$ of the form $\prod_{k<n}I_k$, where
$I_k\in\I_k$ for each $k<n$. For each $k<n$, let $-1\loe
i_0(k)\loe m$ be such that $a^k_{i_0(k)}\loe p_k\loe
a^k_{i_0(k)+1}$, where $a^k_{-1}=0$, $a^k_{m+1}=1$. Define
$r_k=\max\{a^k_{i_0(k)},\inf I_k\}$ and
$s_k=\min\{a^k_{i_0(k)+1},\sup I_k\}$. Pick any $k\in T$. Assume
first that $-1<i_0(k)<m$, so in particular $r_k,s_k\in D_k$ and
hence $r_k,s_k\in\ovr{H}_k$. If $q_1$, $q_2$ denote the
modifications of $p$ obtained by replacing the $k$-th coordinate
of $p$ by $r_k$ and $s_k$ respectively, then
$|g(q_i)-f(q_i)|<(7r+1)\eps$, $i=1,2$, by the induction
hypothesis. Using ($**$) with $j_0=1$ and ``$a^k_{i_0(k)}\loe
r_k\loe p_k\loe s_k\loe a^k_{i_0(k)+1}$'' in the place of
``$a^k_{i_0}=b_0\loe b_1\loe x\loe b_2\loe b_3=a^k_{i_0+1}$'' we
get
\begin{eqnarray}
|g(q_1)-g(p)| + |g(p)-g(q_2)| & \loe & \norm{g(h_k(r_k
,\,\cdot\,))-g(h_k(p_k,\,\cdot\,))}_\infty \nonumber\\
& &+
\norm{g(h_k(p_k,\,\cdot\,))-g(h_k(s_k,\,\cdot\,))}_\infty\nonumber\\
& < & \norm{f(h_k(r_k,\,\cdot\,))-f(h_k(p_k,\,\cdot\,))}_\infty \nonumber\\
& &+\norm{f(h_k(p_k,\,\cdot\,))-f(h_k(s_k,\,\cdot\,))}_\infty +
4\eps\nonumber\\
& < & 6\eps\nonumber
\end{eqnarray}
and hence
\begin{eqnarray}
|g(p)-f(p)| & \loe & |g(p) - g(q_1)| + |g(q_1)-f(q_1)| + |f(q_1)-f(p)|\nonumber\\
& < & 6\eps + (7r+1)\eps + \eps = (7(r+1)+1)\eps.\nonumber
\end{eqnarray}
Assume now that $i_0(k)=m$ (the case $i_0(k)=-1$ is similar). We
have $a^k_0\leq\dots\leq a^k_m\leq p_k$. Let $q$ denote the
modification of $p$ obtained by replacing the $k$-th coordinate
with $a^k_m$. Then
\begin{align*}
|f(q)-f(p)|&\loe\norm{f(h_k(a^k_m,\argum))-f(h_k(p_k,\argum))}_\infty\\
&\loe v^k_{m+1}(f)-\sum_{i<m}\norm{f(h_k(a^k_i,\argum))-
   f(h_k(a^k_{i+1},\argum))}_\infty\\
&<v^k_{m+1}(f)-v^k_m(f)+\delta<2\eps.
\end{align*}
Similarly, using ($*_1$), we get
\begin{align*}
|g(q)-g(p)|&\loe\norm{g(h_k(a^k_m,\argum))-g(h_k(p_k,\argum))}_\infty\\
&\loe
v^k_{m+1}(g)-\sum_{i<m}\norm{g(h_k(a^k_i,\argum))-
   g(h_k(a^k_{i+1},\argum))}_\infty\\
&< v^k_{m+1}(f)+ \eps - \sum_{i<m}|g(x^k_{a^k_i,a^k_{i+1}})-
   g(y^k_{a^k_i,a^k_{i+1}})|\\
&< v^k_{m+1}(f)+ \eps - \sum_{i<m}|f(x^k_{a^k_i,a^k_{i+1}})-
   f(y^k_{a^k_i,a^k_{i+1}})| +\eps\\
&< v^k_{m+1}(f) - \sum_{i<m}\norm{f(h_k(a^k_i,\argum))-
   f(h_k(a^k_{i+1},\argum))}_\infty
+m\delta + \eps\\
&<v^k_{m+1}(f)-v^k_m(f)+(m+1)\delta+2\eps<4\eps.
\end{align*}
Thus $|f(p)-g(p)|<6\eps+|f(q)-g(q)|$ and by the induction
hypothesis, $|f(q)-g(q)|<(7r+1)\eps$. Hence also in this case we
get $|f(p)-g(p)|<(7(r+1)+1)\eps$.

Finally, $\norm{f-g}_\infty<(7n+1)\eps$ which completes the proof.
\end{pf}

\begin{uwgi}
The above result is no longer valid if we drop the requirement
that the sets $D_i$ contain all pairs of adjacent points. For
example, if $L$ is the double arrow line and $D$ is a countable
dense set then $\Tau_p(D)$ is second countable, while $C(L)$ is
not second countable, and the same is true when restricted to any
sphere of $C(L)$.

We also cannot replace the assumption on the sets $D_i$ by ``dense
countably compact''. It is shown in \cite[Example~5.17]{Bu1} that
the space of continuous functions on
$D=(\omega_1+\omega_1^*)^{\omega_1}$ endowed with the topology
induced by the lexicographic order ($\omega_1^*$ means $\omega_1$
with the reversed order) is not norm-SLD for the pointwise
topology. In particular, it has no $\Tau_p$-Kadec renorming. On
the other hand, $D$ is a countably compact linearly ordered space.
If we take $L$ to be the \v{C}ech-Stone compactification of $D$,
then $L$ is linearly ordered---it is obtained from the Dedekind
completion of $D$ by doubling the points which are not endpoints
and are not in $D$---and $C(L)$ is isomorphic to $C(D)$ via the
restriction map. Since this map is also a
$(\Tau_p(D),\Tau_p)$-homeomorphism, $C(L)$ has no
$\Tau_p(D)$-Kadec renorming.
\end{uwgi}

\section{Inverse limits and projectional resolutions of the identity}

In this section we show the existence of a $\tau_p$-Kadec
renorming on a space $C(K)$ when $K$ is a suitable inverse limit
of spaces $K'$ for which $C(K')$ has a $\Tau_p$-Kadec renorming.
As an application, we obtain in particular that $C(K\times L)$ has
a $\Tau_p$-Kadec renorming, whenever $C(K)$ has a $\Tau_p$-Kadec
norm and $L$ is a Valdivia compact space.

We begin with a technical lemma inspired by a very useful result
of Troyanski. (See \cite[VII~Lemma~1.1]{DGZ}.)

\begin{lm}\label{modrzew}
Let $(X,\anorm)$ be a Banach space and let $\Tau$ be a linear
topology on $X$ such that the unit ball of $X$ is
$\Tau$-closed. Fix a function $h\colon \N\to\N$. Suppose there are
\begin{enumerate}
\item[{\rm(a)}]
families $\Ef_0,\Ef_1,\dots$ of bounded $(\Tau,\Tau)$-continuous
linear operators on $X$ such that for each $n$, $\Ef_n$ is
uniformly bounded,
\item[{\rm(b)}]  for each $T\in\bigcup_{n\in\N}\Ef_n$, an equivalent
$\Tau$-Kadec norm $\aabs_T$ on the range of\/ $T$ such that $T$ is
$(\Tau,\Tau)$-continuous, and
\item[{\rm(c)}] for each $n\in\N$ and $T\in \Ef_n$, a set $S_n(T)\subs
\Ef_0\cup\dots\cup\Ef_n$ of cardinality at most $h(n)$,
\end{enumerate}
so that
\begin{enumerate}
\item[{\rm(d)}] for each $x\in X$ and each $\eps>0$, we can find
$n\in\N$ and $T\in \Ef_n$ such that $\norm{x-T_0 x}<\eps$ for some
$T_0\in S_n(T)$ and
$\abs{Tx}_T>\sup\setof{\abs{T'x}_{T'}}{T'\in\Ef_n,\,T'\not=T}$.
\end{enumerate}
Then there exists an equivalent $\Tau$-Kadec norm on $X$.
\end{lm}

\begin{pf}
We may assume that $\aabs_T\loe\anorm$ for each
$T\in\bigcup_{\ntr}\Ef_n$. Define
$$\abs x_{k,n} = \sup\setof{\abs{Tx}_T + { {\frac1{k}}}\sum_{T'\in
S_n(T)}\abs{T'x}_{T'}+\norm{x-T'x}}{T\in \Ef_n}$$ and
$$\abs x =\norm x +\sum_{k,n<\omega}\beta_{k,n}\abs x_{k,n},$$
where $\beta_{k,n}>0$ are such that $\beta_{k,n}\abs
x_{k,n}\loe2^{-(k+n)}\norm x$. (These constants exist because for
each fixed $n$, the operators in $\Ef_n$ are uniformly bounded and
the sets $S_n(T)$, $T\in\Ef_n$, are bounded in cardinality.)

It is clear that $\aabs$ is equivalent to $\|\cdot\|$. We will
show that $\aabs$ is $\Tau$-Kadec. It is $\Tau$-lsc since $\anorm$
and all the $\aabs_{k,n}$ are (use (c) and
Lemma~\ref{alexandrov}). Thus, by Corollary~\ref{snieg1}(b), on
$S:=\setof{x\in X}{\abs x=1}$, each of these functions is
$\Tau$-continuous. Fix $x\in S$ and $\eps>0$. By (d), there are
$n\in \N$ and $T\in \Ef_n$ such that $\norm{x-T_0 x}<\eps$ for
some $T_0\in S_n(T)$ and
$$\delta=\abs{Tx}_T-\sup\setof{\abs{T'x}_{T'}}{T'\in\Ef_n,\,T'\ne T}>0.$$
Choose $k$ so that
$$\frac{h(n)}{k}\sup\setof{2\|T'\|+1}{T'\in\Ef_0\cup\dots\cup\Ef_n}
   \cdot\norm x<\delta.$$
Then
$$
\abs x_{k,n} = \abs{Tx}_T + {\frac1{k}}\sum_{T'\in
S_n(T)}\abs{T'x}_{T'}+\norm{x-T'x}.
$$
(To see this, consider the effect on the expression on the
right-hand side of the equation of replacing $T$ by some other
$\widetilde T\in\Ef_n$. The first term drops by at least $\delta$
(by definition of $\delta$). By the choice of $k$, the second term
cannot make up for the decrease.) By Proposition~\ref{annulus} and
the $(\Tau,\Tau)$-continuity of $T_0$, there is an $\eta>0$ and
there is a $U\in\Tau$ containing $x$ such that if
$\abs{T_0y}_{T_0}$ is within $\eta$ of $\abs{T_0x}_{T_0}$ and
$y\in U$ then $\norm{T_0y-T_0x}<\eps$.

From the $\Tau$-lsc of each of the terms in the expression for
$\abs x_{k,n}$ as functions of $x$ and the $\Tau$-continuity of
$\aabs_{k,n}$ on $S$, it follows from Corollary~\ref{lsc2} that
$y\mapsto\abs{T_0y}_{T_0}$ and $y\mapsto\norm{y-T_0y}$ are
continuous at $x$ on $S$. Thus, by shrinking $U$ to a smaller
$\Tau$-neighborhood of $x$, we may arrange that
$y\mapsto\abs{T_0y}_{T_0}$ and $y\mapsto\norm{y-T_0y}$ vary by
less than $\min\{\eta,\eps\}$ on $U\cap S$. Since
$\|x-T_0x\|<\eps$, this means in particular that
$\|y-T_0y\|<2\eps$ for $y\in U\cap S$.

For $y\in U\cap S$, we have
$$
\norm{y-x}\loe\norm{y-T_0y}+\norm{T_0y-T_0x}+\norm{T_0x-x}<
   2\eps+\eps+\eps=4\eps.
$$
This completes the proof.
\end{pf}

\begin{uwgi}
The above lemma, as well as its corollaries, could be stated in a
more general form saying that on each $T X$ there is a weaker
linear topology $\Tau_T$ for which $T$ is
$(\Tau,\Tau_T)$-continuous and $T X$ has a $\Tau_T$-Kadec
renorming. The proofs require only minor changes.
\end{uwgi}

\begin{tw}
Let $X$ be a Banach space and let $\sett{\map{P_n}XX}{n\in\N}$ be
a uniformly bounded sequence of projections such that
$\bigcup_{n\in\N}P_nX$ is dense in $X$. Let $\Tau$ be a weaker
linear topology on $X$ such that the unit ball is $\Tau$-closed.
If for each $n\in\N$, $P_n$ is $(\Tau,\Tau)$-continuous and there
exists a $\Tau$-Kadec renorming of $P_nX$, then there exists a
$\Tau$-Kadec renorming of $X$.
\end{tw}

\begin{pf}
We apply Lemma \ref{modrzew} with $\Ef_n=\sn{P_n}$ and
$S_n(P_n)=\sn{P_n}$. Condition (d) of Lemma \ref{modrzew} reduces
in this case to the fact that for every $x\in X$ and $\eps>0$
there exists $\ntr$ such that $\norm{x-P_nx}<\eps$. To see that
this is true, fix $x\in X$ and $\eps>0$ and set
$\delta=\eps/(1+M)$, where $M$ is a constant which bounds the
norms of all $P_n$'s. Then, by assumption, there are $n\in\N$ and
$y\in P_nX$ such that $\norm{x-y}<\delta$. We have $y=P_ny$ and
hence $\norm{P_nx-y}\loe\norm{P_n}\cdot\norm{x-y}<M\delta$. Thus
$$\norm{x-P_nx}\loe \norm{x-y}+\norm{y-P_nx}<\delta+M\delta=\eps.$$
\end{pf}

\begin{tw}\label{jesion}
Let $(X,\anorm)$ be a Banach space and assume that
$\sett{\map{T_\al}XX}{\al<\kappa}$ is a sequence of uniformly
bounded linear operators on $X$ such that for each $x\in X$,
\begin{enumerate}
    \item[{\rm(i)}] the sequence $\sett{\|T_\al x\|}{\al<\kappa}$
               belongs to $c_0(\kappa)$,
    \item[{\rm(ii)}] for every $\eps>0$ there exists a finite set
                $A\subs\kappa$ such that
                $$\bnorm{x-\sum_{\al\in A}T_\al x}<\eps,$$
    \item[{\rm(iii)}] $T_\al X\cap T_\beta X=\sn 0$ whenever
                 $\al\ne \beta$.
\end{enumerate}
Assume further that $\Tau$ is a linear topology on $X$ such that
the the unit ball of $X$ is $\Tau$-closed and for each
$\al<\kappa$, $T_\al X$ has a $\Tau$-Kadec renorming and $T_\al$
is $(\Tau,\Tau)$-continuous. Then $X$ has an equivalent
$\Tau$-Kadec norm.
\end{tw}

\begin{pf}
Let $Q_A=\sum_{\al\in A}T_\al$ and define
$\Ef_n=\setof{Q_A}{A\in\dpower\kappa n}$,
$S(Q_A)=\setof{Q_{A'}}{A'\subs A}$ (so  $S(Q_A)$ has cardinality
at most $2^{|A|}$).
If $\anorm_\al$ is a $\Tau$-Kadec norm on $T_\al X$ then
$\anorm_{Q_A}=\sum_{\al\in A}\anorm_\al$ is a $\Tau$-Kadec
norm on $Q_AX$. We may assume that $\anorm_\al\loe\anorm$ for each
$\al<\kappa$. We need to check condition (d) of Lemma
\ref{modrzew}. Fix $x\in X$, $\eps>0$. By (ii) there exists
$A_0\in\fin\kappa$ such that $\norm{x-Q_{A_0}x}<\eps$. By (i),
there exists a finite set $A\sups A_0$ such that
$$\max_{\al\notin A}\norm{T_\al x}_\al<\min_{\al\in A}\norm{T_\al x}_\al.$$
It follows that
$\norm{Q_Ax}_{Q_A}>\sup\setof{\norm{Q_Bx}_{Q_B}}{|B|=|A|\Land B\ne
A}$. Thus, by Lemma \ref{modrzew}, we get a $\Tau$-Kadec renorming
of $X$.
\end{pf}

\begin{tw}\label{goryl}
Assume $X$ is a Banach space and
$\sett{\map{P_\al}XX}{\al\loe\kappa}$ is a sequence of projections
such that
\begin{enumerate}
    \item[{\rm(a)}] $P_0=0$, $P_\kappa=\id_E$ and $P_\beta P_\al=P_\al=P_\al P_\beta$
          whenever $\al\loe\beta\loe\kappa$.
    \item[{\rm(b)}] There is $M<+\infty$ such that $\norm{P_\al}\loe M$ for every
          $\al<\kappa$.
    \item[{\rm(c)}] If $\lam\loe\kappa$ is a limit ordinal then
          $\bigcup_{\xi<\lam}P_\xi E$ is dense in $P_\lam E$.
\end{enumerate}
Assume that $\Tau$ is a linear topology on $X$ such that the unit
ball of $X$ is $\Tau$-closed and for each $\al<\kappa$,
$(P_{\al+1}-P_\al)X$ has a $\Tau$-Kadec renorming and
$P_{\al+1}-P_\al$ is $(\Tau,\Tau)$-continuous. Then $X$ has a
$\Tau$-Kadec renorming.
\end{tw}

\begin{pf}
Let $T_\al=P_{\al+1}-P_\al$. A standard and well known argument
(see e.g. \cite[pp. 236, 284]{DGZ}) shows that
$\sett{T_\al}{\al<\kappa}$ satisfies the assumptions of Theorem
\ref{jesion}. We write out the proof of condition (ii) for the
sake of completeness because it is not given explicitly in
\cite{DGZ}.

Proceed by induction on limit ordinals $\lam<\kappa$. If
$\lam=\omega$ then $P_\omega x=\lim_{n\to\infty}P_n
x=\sum_{\ntr}(P_{n+1}x-P_n x)=\sum_{\ntr}T_n x$ (recall that
$P_0=0$), so $\sum_{n<k}T_n x$ can be taken arbitrarily close to
$P_\omega x$. Now let $\lam>\omega$ and assume the statement is
true for limit ordinals below $\lam$ (and for every $\eps>0$).
There exists $\xi_0<\lam$ such that $\norm{P_\lam x- P_\beta
x}<\eps/2$ for $\xi\goe\xi_0$. If there is a limit ordinal $\beta$
such that $\xi_0\loe\beta<\lam$ then, by induction hypothesis
\begin{equation}
\bnorm{P_\beta x-\sum_{\al\in A}T_\al x}<\eps/2.\tag{$*$}
\end{equation}
for some finite set $A\subs\beta$ and we have $\norm{P_\lam
x-\sum_{\al\in A}T_\al x}<\eps$. Otherwise, $\xi_0=\beta+n$, where
$\beta\goe\omega$ is a limit ordinal and again ($*$) holds for
some finite set $A\subs\beta$. Now we have $P_{\beta+n}x-P_\beta
x=\sum_{\al=\beta}^{\beta+n-1}T_\al x$ and hence $\norm{P_\lam
x-\sum_{\al\in B}T_\al x}\loe \norm{P_\lam
x-P_{\beta+n}x}+\norm{P_\beta x-\sum_{\al\in A} T_\al x}<\eps$,
where $B=A\cup\{\beta,\beta+1,\dots,\beta+n-1\}$.
\end{pf}

A sequence $\setof{P_{\al}}{\al\loe\kappa}$ satisfying conditions
(a), (b) and (c) of the above theorem with $M=1$ and such that the
density of $P_\al X$ is $\loe|\al|+\aleph_0$, is called a {\em
projectional resolution of the identity} (PRI) on $X$, see
\cite{DGZ} or \cite{Fa}.

The following proposition is a purely category-theoretic property
of inverse limits.  It is standard but we do not know a reference
for it, so we write out the proof.

\begin{prop}\label{retrakcje}
Let $\setof{X_\al; p^\beta_\al}{\al<\beta<\kappa}$ be a continuous
inverse sequence of topological spaces such that each
$p^{\al+1}_\al$ is a retraction and let $X$, with projections
$\setof{p_\al}{\al<\kappa}$, be the inverse limit of the sequence.
Then there exists a collection of continuous embeddings
$\sett{\map{i^\beta_\al}{X_\al}{X_\beta}}{\al<\beta<\kappa}$, such
that
\begin{enumerate}
  \item[{\rm(1)}]
$p^\beta_\al i^\beta_\al = \id_{X_\al}$ for all $\al<\beta<\kappa$
and $\;i^\beta_\gamma i^\gamma_\al = i^\beta_\al$ for all
$\al<\gamma<\beta<\kappa$.
\end{enumerate}
Moreover, there exist continuous embeddings $\map{i_\al}{X_\al}X$
such that
\begin{enumerate}
    \item[{\rm(2)}]
$p_\al i_\al = \id_{X_\al}$ and $\;i_\beta i^\beta_\al = i_\al$,
whenever $\al<\beta<\kappa$.
\end{enumerate}
\end{prop}

\begin{pf}
We can treat (2) as a special case of (1) by allowing
$\beta=\kappa$ in (1) and setting $X_\kappa=X$ and
$p^\kappa_\al=p_\al$ for $\al<\kappa$. We construct the maps
$i_\al^\beta$ by induction on $\beta\leq\kappa$. Assume
$i_\xi^\eta$ have been constructed for every $\xi<\eta<\beta$; for
convenience we set $i^\xi_\xi=\id_{X_\xi}$. Suppose first that
$\beta$ is a successor, i.e. $\beta=\delta+1$.

Fix any continuous map $\map{
i^{\delta+1}_\delta}{X_\delta}{X_{\delta+1}}$ which is a right
inverse of $p^{\delta+1}_\delta$. For $\al<\delta$, define
$i_\al^{\delta+1} = i^{\delta+1}_\delta i_\al^\delta$. To see that
(1) holds, observe that
$$
p^{{\delta+1}}_\al i^{\delta+1}_\al = p^{\delta+1}_\al
i^{\delta+1}_\delta i^\delta_\al = p_\al^\delta
p_\delta^{\delta+1} i^{\delta+1}_\delta i_\al^\delta =
p_\al^\delta i_\al^\delta=\id_{X_\al},
$$
and
$$
i^{\delta+1}_\gamma i^\gamma_\al = i^{\delta+1}_\delta
i^\delta_\gamma i^\gamma_\al = i^{\delta+1}_\delta i^\delta_\al =
i^{\delta+1}_\al.
$$

Suppose now that $\beta$ is a limit ordinal. Fix $\al<\beta$.
Observe that for $\al\loe\xi<\eta<\beta$ we have
$$
p^\eta_\xi i^\eta_\al = p^\eta_\xi i^\eta_\xi i^\xi_\al =
i^\xi_\al.
$$
Since $X_\beta$ together with
$\setof{p^\beta_\xi}{\xi\in[\al,\beta)}$ is the limit of
$\setof{X_\xi; p^\eta_\xi}{\al\loe\xi<\eta<\beta}$,
%($[\al,\beta)$
%denotes the interval of ordinals smaller than $\beta$ and not
%smaller than $\al$),
there exists a unique continuous map
$\map{i^\beta_\al}{X_\al}{X_\beta}$ such that
$$
p^\beta_\xi i^\beta_\al = i^\xi_\al
$$
holds for every $\xi\in[\al,\beta)$. In particular $p^\beta_\al
i^\beta_\al=\id_{X_\al}$. Thus we have defined mappings
$i^\beta_\al$, for $\al<\beta$. It remains to check that
$i^\beta_\gamma i^\gamma_\al = i^\beta_\al$ for
$\al<\gamma<\beta$. To see this, observe that for
$\xi\in[\gamma,\beta)$ we have
$$
p^\beta_\xi (i^\beta_\gamma i^\gamma_\al) = i^\xi_\gamma
i^\gamma_\al = i^\xi_\al,
$$
and for $\xi\in[\al,\gamma)$ we have
$$
p^\beta_\xi (i^\beta_\gamma i^\gamma_\al) = p^\gamma_\xi
p^\beta_\gamma (i^\beta_\gamma i^\gamma_\al) = p^\gamma_\xi
i^\gamma_\al = i^\xi_\al.
$$
Since $i^\beta_\al$ is the unique map satisfying $p^\beta_\xi
i^\beta_\al = i^\xi_\al$ for $\xi\in[\al,\beta)$, we get
$i^\beta_\gamma i^\gamma_\al=i^\beta_\al$.
\end{pf}

\begin{lm}\label{kot}
Let $\{K;p_\al:\al<\kappa\}$ be the inverse limit of the
continuous inverse sequence of compact spaces
$$\setof{K_\al;p^{\beta}_\al}{\al<\beta<\kappa}$$
in which the bonding maps $p^{\al+1}_\al$ are retractions.
\begin{enumerate}
\item[{\rm(a)}]
If for each $\al<\kappa$, $C(K_\al)$ has a $\Tau_p$-Kadec
renorming, then $C(K)$ has a $\Tau_p$-Kadec renorming.
\item[{\rm(b)}]
Let $\setof{i^\beta_\al}{\al<\beta<\kappa}$ and
$\setof{i_\al}{\al<\kappa}$ be collections of right inverses
satisfying {\rm(1)} and {\rm(2)} of Proposition~\ref{retrakcje}.
Assume that $D\subs K$ is dense, and for each $\al<\kappa$, $i_\al
p_\al[D]\subseteq D$ and $C(K_\al)$ has a $\Tau_p(p_\al[D])$-Kadec
renorming. Then $C(K)$ has a $\Tau_p(D)$-Kadec renorming.
\end{enumerate}
\end{lm}

\begin{pf}
(a) (Cf.\ the proof of \cite[VI~Theorem~7.6]{DGZ}.) Let
$\setof{i^\beta_\al}{\al<\beta<\kappa}$ and
$\setof{i_\al}{\al<\kappa}$ be collections of right inverses given
by Lemma \ref{retrakcje}. Let $R_\al=i_\al p_\al$. $R_\al$ is a
retraction of $K$ onto $i_\al[K_\al]$. If $\al<\beta$ then
$$
R_\al R_\beta = i_\al\, p_\al\, i_\beta\, p_\beta = i_\al\,
p^\beta_\al\, p_\beta\, i_\beta\, p_\beta = i_\al\, p^\beta_\al\,
p_\beta = i_\al\, p_\al = R_\al
$$
and
$$
R_\beta R_\al = i_\beta\, p_\beta\, i_\al\, p_\al = i_\beta\,
p_\beta\, i_\beta\, i^\beta_\al\, p_\al = i_\beta\, i^\beta_\al\,
p_\al = i_\al\, p_\al = R_\al.
$$
We also have $R_\al R_\al=R_\al$.

Let $\map{P_\al}{C(K)}{C(K)}$ be given by $P_\al(f)=f R_\al$.

$C(K_\al)$ can be identified with the range of $P_\al$ via the
linear map $T$ defined by $Tg = gp_\al$. $T$ is norm-preserving,
and in particular one-to-one, because $p_\al$ maps onto $K_\al$.
From $ Tg = gp_\al=gp_\al i_\al p_\al=gp_\al R_\al =
P_\al(gp_\al)$ and $P_\al(f)=fR_\al = f i_\al p_\al = T(f i_\al)$,
we see that the range of $T$ is indeed the same as the range of
$P_\al$. Note that $T^{-1}(h)=hi_\al$. $T$ is a
$\Tau_p$-homeomorhism because for $x\in K$ and $y\in K_\al$, the
maps $g\mapsto (Tg)(x)=g(p_\al x)$ and $h\mapsto
(T^{-1}h)(y)=h(i_\al y)$ are $\Tau_p$-continuous. It follows from
our assumption that the range of $P_\al$ has an equivalent
$\Tau_p$-Kadec norm.

Then $\setof{P_\al}{\al<\kappa}$ is a sequence of projections of
norm one satisfying the condition
$$
\al<\beta \implies P_\al P_\beta = P_\beta P_\al = P_\al.
$$
For any $x\in K$, the map $f\mapsto P_\al(f)(x)$ is
$\Tau_p$-continuous since it coincides with $x\mapsto f(R_\al x)$.
Hence, $P_\al$ is $(\Tau_p,\Tau_p)$-continuous.

We now check that $\bigcup_{\al<\beta}P_\al C(K)$ is dense in
$P_\beta C(K)$ for every limit ordinal $\beta\loe \kappa$. It will
then follow that $\sett{P_\al}{\al<\kappa}$ satisfies the
assumptions of Theorem \ref{goryl} and the proof of (a) will be
complete. We show that for each $f\in C(K)$,
\[
\lim_{\al\to\beta}P_\al(f)=P_\beta(f).
\]
Fix $\eps>0$. $K_\beta$ has a base consisting of open sets of the
form $(p^\beta_\al)^{-1}[U]$ where $\al<\beta$ and $U$ is open in
$K_\al$. Hence, $K_\beta$ is covered by finitely many such sets on
which the oscillation of $f\, i_\beta$ is at most $\eps$. By
replacing the finitely many $\al$'s involved here by the largest
of them, we may assume that they are all equal to some
$\al_0<\beta$. (If $\al<\al'<\beta$ and $U$ is open in $K_\al$,
then since $p^\beta_\al = p^{\al'}_\al p^{\beta}_{\al'}$, we have
$(p^\beta_\al)^{-1}[U]=(p^\beta_{\al'})^{-1}[V]$ where
$V=(p^{\al'}_\al)^{-1}[U]$.) Thus we have open sets
$U_1,\dots,U_n$ in $K_{\al_0}$ such that the sets
\[
(p^\beta_{\al_0})^{-1}[U_1],\dots,(p^\beta_{\al_0})^{-1}[U_n]
\]
cover $K_\beta$ and on each of them the oscillation of $f\,
i_\beta$ is at most $\eps$. For any $\al$ such that $\al_0\leq \al
<\beta$ and for any $x\in K$, letting $j\in\{1,\dots,n\}$ be such
that $p_{\al_0}(x)\in U_j$, we have
\[
p^\beta_{\al_0}(i^\beta_\al p_\al(x)) = p^\al_{\al_0}\, p_\al(x) =
p_{\al_0}(x)\in U_j,
\]
so that
\[
i^\beta_\al\, p_\al(x)\in (p^\beta_{\al_0})^{-1}[U_j].
\]
Clearly we also have
\[
p_\beta(x)\in (p^\beta_{\al_0})^{-1}[U_j],
\]
and hence
\[
|P_\al(f)(x)-P_\beta(f)(x)|=|fR_\al(x)-fR_\beta(x)|= |f\, i_\al\,
p_\al (x) - f\, i_\beta\, p_\beta (x)| = |f\, i_\beta (i^\beta_\al
p_\al (x)) - f\, i_\beta (p_\beta (x))| \leq \eps.
\]
This completes the proof of (a).

The proof of (b) is obtained by making suitable adjustments to the
proof of (a). We check that $T$ is a
$(\Tau_p(p_\al[D]),\Tau_p(i_\al\,p_\al[D]))$-homeomorphism. When
$d\in D$, $g\mapsto
(Tg)(i_\al\,p_\al(d))=g(p_\al\,i_\al\,p_\al(d)) =g(p_\al(d))$ is
$\Tau_p(p_\al[D])$-continuous and $h\mapsto
(T^{-1}h)(p_\al(d))=h(i_\al\,p_\al(d))$ is
$\Tau_p(i_\al\,p_\al[D])$-continuous. Hence, our assumption gives
that the range of $P_\al$ has an equivalent
$\Tau_p(i_\al\,p_\al[D])$-Kadec norm. It follows that the range of
$P_\al$ has an equivalent $\Tau_p(D)$-Kadec norm. For any $d\in
D$, the map $f\mapsto P_\al(f)(d)=f(R_\al(d))=f(i_\al\,p_\al(d))$
is $\Tau_p(D)$-continuous. Hence, $P_\al$ is
$(\Tau_p(D),\Tau_p(D))$-continuous. Finally, the fact that $D$ is
dense ensures that the unit ball of $C(K)$ is $\Tau_p(D)$-closed.
The rest of the proof is as for (a).
\end{pf}

Given a family of spaces $\sett{X_\al}{\al<\kappa}$, their product
$\prod_{\al<\kappa}X_\al$ is the limit of a continuous inverse
sequence of smaller products $\prod_{\xi<\al}X_\xi$, with the
usual projections as bonding maps. This leads to the following.

\begin{wn}
Let $\setof{K_\al}{\al<\kappa}$ be a family of compacta and assume
that for every finite $S\subs\kappa$, $C(\prod_{\al\in S}K_\al)$
has a $\Tau_p$-Kadec renorming. Then $C(\prod_{\al<\kappa}K_\al)$
has a $\Tau_p$-Kadec renorming.
\end{wn}

\begin{pf}
Proceed by induction on the cardinality of the index set, which we
can assume is infinite. The induction hypothesis ensures that for
each $\beta<\kappa$, $C(\prod_{\al<\beta}K_\al)$ has a
$\Tau_p$-Kadec renorming. Now apply Lemma~\ref{kot}(a).
\end{pf}

In \cite{JNR3} an analogous result on the $\sig$-fragmentability
of products is proved. In \cite{BR} it is shown that the property
of having a $\Tau_p$-lsc LUR renorming is productive in the sense
that $C(\prod_{\al<\kappa}K_\al)$ has a $\Tau_p$-lsc LUR renorming
if (and trivially only if) each $C(K_\al)$ has a $\Tau_p$-lsc LUR
renorming. It is unknown whether the property of having a
$\Tau_p$-Kadec renorming is productive in this sense.

Lemma \ref{kot} allows us to generalize Theorem \ref{wiatr} to
infinite products.

\begin{tw}\label{products of lo}
Let $\setof{L_\al}{\al<\kappa}$ be a collection of compact
linearly ordered spaces and for each $\al<\kappa$ let $D_\al$ be a
dense subset of $L_\al$ which contains all pairs of adjacent
points. Then $C(\prod_{\al<\kappa}L_\al)$ has an equivalent
$\Tau_p(\prod_{\al<\kappa}D_\al)$-Kadec norm.
\end{tw}

\begin{pf}
Proceed by induction on the cardinality of the index set.
Theorem~\ref{wiatr} takes care of the case $\kappa<\omega$. Assume
that $\kappa$ is an infinite cardinal and write
$K=\prod_{\al<\kappa}L_\al$ and $K_\al=\prod_{\xi<\al}L_\xi$ for
$\al<\kappa$. Note that $K$, equipped with the usual projections
$p_\al\colon K\to K_\al$, is the inverse limit of the continuous
inverse sequence $\setof{K_\al;p^\beta_\al}{\al<\beta<\kappa}$,
where the $p^\beta_\al$'s are the usual projections. Fix a base
point $d_\al\in D_\al$ for each $\al<\kappa$. For
$\al<\beta<\kappa$, define embeddings
\[
\textstyle{i^\beta_\al\colon \prod_{\xi<\al}L_\xi\to
\prod_{\xi<\beta}L_\xi}
\]
by $i^\beta_\al(x)(\xi)=x(\xi)$ for $\xi<\al$ and
$i^\beta_\al(x)(\xi)=d_\xi$ for $\al\leq \xi<\beta$. By the
induction hypothesis, $C(K_\al)$ has an equivalent
$\Tau_p(\prod_{\xi<\al}D_\xi)$-Kadec norm for each $\al,\kappa$.
The assumptions of part (b) of Lemma~\ref{kot} are satisfied with
$D=\prod_{\al<\kappa}D_\al$.
\end{pf}

Denote by $\R$ the minimal class of compact spaces which contains
all metric compacta and is closed under limits of continuous
inverse sequences of retractions. More formally, $\R$ is the
smallest class of spaces which satisfies the following conditions:
\begin{enumerate}
    \item Every metrizable compact space is in \R.
    \item If $\S=\setof{X_\al; p^\beta_\al}{\al<\beta<\kappa}$
          is a continuous inverse sequence such that each $X_\al$ is in
          \R~and each $p^{\al+1}_\al$ is a retraction, then every
          space homeomorphic to $\liminv\S$ belongs to $\R$.
\end{enumerate}
Note that every Valdivia compact space belongs to \R~ (see e.g.
\cite{Ka}). Also, for every ordinal $\xi$, the compact linearly
ordered space $\xi+1$ belongs to \R. If $\xi\goe\aleph_2$ then
$\xi+1$ is not Valdivia compact (see \cite{Ka}). It is easy to see
that class \R~is closed under products and direct sums.

\begin{tw}\label{product inv lim}
{\rm(a)} Assume $K$ is a compact space such that $C(K)$ has a
$\Tau_p$-Kadec renorming and assume $L\in\R$. Then $C(K\times L)$
has a $\Tau_p$-Kadec renorming.

{\rm(b)} For every $L\in\R$, $C(L)$ has a $\Tau_p$-lsc LUR
renorming.
\end{tw}

\begin{pf}
(a) Denote by $\R_0$ the class of all spaces $L\in\R$ such that
$C(K\times L)$ has a $\Tau_p$-Kadec renorming. It suffices to show
that $\R_0$ contains all metric compacta and is closed under
limits of continuous inverse sequences of retractions. The latter
fact follows from Lemma \ref{kot}, because if
$L=\liminv\setof{L_\al;p^\beta_\al}{\al<\beta<\kappa}$ then
$K\times L=\liminv\setof{K\times L_\al;
q^\beta_\al}{\al<\beta<\kappa}$, where $q^\beta_\al=\id_K\times
p^\beta_\al$. It remains to show that $\R_0$ contains all metric
compacta. As every compact metric space is a continuous image of
the Cantor set, it is enough to show that $C(K\times 2^\omega)$
has a $\Tau_p$-Kadec renorming.

We have $2^\omega=\liminv\setof{2^n;p^m_n}{n<m<\omega}$ so
$$K\times 2^\omega=\liminv\setof{K\times 2^n;q^m_n}{n<m<\omega},$$
where $q^m_n=\id_K\times p^m_n$. Clearly, $C(K\times 2^n)$ has a
$\Tau_p$-Kadec renorming being a finite power of $C(K)$, so again
Lemma \ref{kot} gives a $\Tau_p$-Kadec renorming of
$C(K\times2^\omega)$.

(b) It is enough to check that the class of all compact spaces
$K$ for which $C(K)$ has a $\Tau_p$-lsc LUR renorming is closed
under inverse limits of retractions. Assume $K=\liminv\S$, where
$\S=\setof{K_\al; r_\al^\beta}{\al<\beta<\kappa}$ is a continuous
inverse sequence of retractions and for each $\al<\kappa$,
$C(K_\al)$ has a $\tau_p$-lsc LUR renorming. As in the proof of
Lemma \ref{kot}(a), there is a sequence of projections
$\setof{P_\al}{\al<\kappa}$ on $C(K)$ such that $P_\al$ is
adjoint to the retraction $\map{r_\al}{K}{K_\al}$. Now apply
Proposition VII.1.6 and Remark VII.1.7 from \cite{DGZ} to obtain
a $\tau_p$-lsc LUR renorming of $C(K)$. In fact,
\cite[Proposition VII.1.6]{DGZ} deals with projectional
resolutions of the identity, but no assumption about the density
of $\im P_\al$ is used in the proof.
\end{pf}

\begin{uwgi}
Note that by Proposition~\ref{continuous image},
Theorem~\ref{product inv lim}(a) applies also when $L$ is a
continuous image of a space from $\R$. (If $L'$ is a continuous
image of $L$, then $K\times L'$ is a continuous image of $K\times
L$.)
\end{uwgi}

\begin{ex}
In \cite{To} an example of a compact, non-separable ccc space of
countable $\pi$-character which has a continuous map onto the
Cantor set in such a way that the fibers are relatively small
linearly ordered spaces (their order type is an ordinal less than
the additivity of Lebesgue measure). This space belongs to \R.
\end{ex}

As in \cite{To}, we use Boolean algebraic language and work with
the Boolean algebra whose Stone space is the required example.

Let \N~denote the set of positive natural numbers and denote by
$\N[i]$ the set of all numbers of the form $2^i(2j-1)$, where
$j\in\N$. Define $K=\setof{x\subs\N}{(\forall\;i)\;|x[i]|\loe i}$,
where $x[i]=x\cap\N[i]$, and
$$Z=\setof{x\in K}{\lim_{i\to\infty}{|x[i]|}/{i}=0}.$$

Denote by $\subs^*$ the {\em almost inclusion relation}, i.e.
$a\subs^* b$ if $a\setminus b$ is finite. Define
$$T=\setof{(t,n)}{n\in\N,\ t\in K\ \mbox{and}\ t\subs n}.$$
We are going to define a subalgebra of $\Pee(T)/\Fin$, where fin
is the ideal of finite subsets of $T$. Let
$$T_{(t,n)}=\setof{(s,m)\in T}{m\goe n\ \mbox{and}\ s\cap n=t}$$
and
$$T_a=\setof{(s,m)\in T}{a\cap m\subs s}.$$
Define $\Be_0$ to be the subalgebra of $\Pee(T)/\Fin$ generated by
the classes of the sets $T_{(t,n)}$, $(t,n)\in T$. Then $\Be_0$ is
a countable free Boolean algebra. In what follows we shall
identify subsets of $T$ with their equivalence classes in
$\Pee(T)/\Fin$. The context should make it clear when classes are
intended.

By \cite[p. 151]{FK}, there exists a sequence
$A=\setof{a_\al}{\al<\kappa}$ of elements of $Z$ such that
$\al<\beta\implies a_\al\subs^* a_\beta$ and for every $a\in K$
there is $\al<\kappa$ such that $a_\al\not\subs^* a$. Moreover
$\kappa$ equals the additivity of the Lebesgue measure, so
$\kappa>\aleph_0$. Let $\Be_\al$ be the subalgebra of
$\Pee(T)/\Fin$ generated by
$$\Be_0\cup\setof{T_a}{a\in K\Land (\exists\; \xi<\al)\; a=^* a_\xi}.$$
Finally, let $\Be=\bigcup_{\al<\kappa}\Be_\al$ and let $X$ be the
Stone space of $\Be$. It has been shown in \cite{To} that $X$ is a
non-separable ccc space with countable $\pi$-character. Moreover,
the inclusion $\Be_0\subs\Be$ induces, by duality, a map from $X$
onto the Cantor set such that all fibers are well-ordered of size
$<\kappa$.

\begin{tw}
$X\in\R$ and consequently $C(X)$ has a $\Tau_p$-lsc LUR renorming.
\end{tw}

\begin{pf}
We will show by induction on $\al<\kappa$ that
$\ult{\Be_\al}\in\R$ for every $\al<\kappa$ and that each quotient
mapping $\map{r_\al}{\ult{\Be_{\al+1}}}{\ult{\Be_\al}}$ induced by
$\Be_\al\subs \Be_{\al+1}$ is a retraction. The latter property is
equivalent to the existence of a retraction $\map
h{\Be_{\al+1}}{\Be_\al}$, i.e. a homomorphism such that $h\rest
\Be_\al=\id_{\Be_\al}$.

Fix $\al<\kappa$ and assume $\ult{\Be_\al}\in\R$. Given Boolean
algebras $\Aa\subs\Be$ and $x\in\Be\setminus \Aa$ we will denote
by $\Aa[x]$ the algebra generated by $\Aa\cup \sn x$ ($\Aa[x]$ is
called a {\em simple extension} of $\Aa$). Note the following

\begin{claim}\label{tekno0}
Assume $a\subs a'$ are in $K$ and $a'\setminus a$ is finite. Then
$T_{a'}\in\Be_0[T_a]$.
\end{claim}

\begin{pf}
Let $\ntr$ be such that $a'\setminus a\subs n$. Let
$\Es=\setof{s\subs n}{s\in K$ and $a'\cap n\subs s}$. Then
$T_{a'}=T_a\cap \bigcup_{s\in\Es}T_{(s,n)}$.
\end{pf}

Define $\Be_{\al+1}^{-1}=\Be_\al$ and
$\Be_{\al+1}^{n+1}=\Be_{\al+1}^n[T_{a_\al\setminus n}]$. By the
above claim,  $\Be_{\al+1}=\bigcup_{\ntr}\Be_{\al+1}^n$. We need
to check that $\Be_{\al+1}^n$ is a retract of $\Be_{\al+1}^{n+1}$
and that $\ult{\Be_\al^{n+1}}\in\R$ for every $n\goe-1$.

Note that, by Sikorski's extension criterion (see e.g.
\cite[p.~67]{Ko}), if $\Aa$ is a Boolean algebra and $\Aa[x]$ is a
simple extension of $\Aa$ then $\Aa$ is a retract of $\Aa[x]$ iff
there exists $c\in \Aa$ such that for every $a_0,a_1\in\Aa$ with
$a_0\loe x\loe a_1$ we have $a_0\loe c \loe a_1$. This holds for
example, if $\setof{a\in\Aa}{a\loe x}$ has a least upper bound in
$\Aa$.

We will need the following easy fact about our Boolean algebra. We
leave the verification to the reader. Part~(a) is like Claim~1
from the proof of \cite[Theorem~8.4]{To}.

\begin{claim}\label{tekno}
{\rm(a)} The sets $T_a\cap T_{(t,n)}$, where $a=^*a_\xi$ for some
$\xi<\al$, are dense in $\Be_\al$.

{\rm(b)} For every nonnegative integer $n$, every element of
$\Be_{\al+1}^n$ is a finite sum of elements of the form $T_a\cap
T_{(t,n)}\cap \neg T_{b_0}\cap\dots\cap\neg T_{b_{k-1}}$, where
$b_i=^* a_{\eta_i}$ for some $\eta_i\loe\al$ and $a=^*a_\xi$ for
some $\xi<\al$ or $a=a_\al\setminus i$ where $i<n$.

{\rm(c)} If $x=T_a\cap T_{(t,n)}$ and $0_{\Be}<x\loe T_{b}$ then
$b\subs a\cup t$.
\end{claim}

We consider separately the cases $n=-1$ and $n>-1$.

{\it Case 1.} $n=-1$. By Claim~\ref{tekno} (a) and (c), no
non-zero element of $\Be_\al$ is below $T_{a_\al}$. Thus
$\Be_\al=\Be_{\al+1}^{-1}$ is a retract of $\Be_{\al+1}^0$. To see
that $\ult{\Be_{\al+1}^0}\in\R$ it is enough to show that
$\Be_{\al}/\Ai$ is countable (and hence its Stone space is second
countable), where $\Ai=\setof{x\in\Be_\al}{x\cap
T_{a_\al}=0_\Be}$, because $\ult{\Be_{\al+1}^0}$ is the direct sum
of $\ult{\Be_{\al}}$ and $\ult{\Be_\al/\Ai}$. Let $\map
q{\Be_{\al}}{\Be_{\al}/\Ai}$ be the quotient map. Observe that for
$\xi<\al$, $q(T_{a_\xi\cap a_\al\setminus n})=1_{\Be_{\al}/\Ai}$,
because $T_{a_\al}\loe T_{a_\xi\cap a_\al\setminus n}$. Now, by
Claim \ref{tekno0}, $\Be_\al$ is generated by
$\Be_0\cup\setof{T_a}{a=a_\xi\cap a_\al\setminus n\Land \ntr\Land
\xi<\al}$. It follows that $\Be_\al/\Ai$ is countable.

{\it Case 2.} $n>-1$. By Claim \ref{tekno}, we have
$\sup\setof{x\in\Be_{\al+1}^n}{x\loe T_{a_\al\setminus
n}}=T_{a_\al\setminus(n-1)}\in\Be_{\al+1}^n$. Hence
$\Be_{\al+1}^n$ is a retract of $\Be_{\al+1}^{n+1}$. In order to
see that $\ult{\Be_{\al+1}^{n+1}}\in\R$ it is enough to show that,
as in Case 1, the quotient algebra $\Be_{\al+1}^n/\Ai$ is
countable, where $\Ai=\setof{x\in\Be_{\al+1}^n}{x\cap
T_{a_\al\setminus n}=0_\Be}$. This can be done by an argument
similar to the one used as in Case 1. We now have new generators
of the form $T_{a_\al\setminus i}$, $i<n$, but only finitely many
of them, so the quotient $\Be_{\al+1}^n/\Ai$ is still countable.
\end{pf}

\begin{uwgi}
If the additivity of the Lebesgue measure is $>\aleph_2$ then the
space $X$ from the above example is not a continuous image of a
Valdivia compact space. Indeed, let $\kappa$ denote the
additivity of the Lebesgue measure and suppose that $X$ is a
continuous image of a Valdivia compact space. Let $\map
hX{2^\omega}$ be a continuous map such that all fibers of $h$ are
well ordered of order type $<\kappa$ (see \cite{To}). One can
show that in fact there are fibers of arbitrary large order type
below $\kappa$ (see the proof of Claim 4 in \cite[p. 74]{To}).
Hence, assuming $\kappa>\aleph_2$, there is $p\in2^\omega$ such
that $F=h^{-1}(p)$ has order type $>\aleph_2$. Observe that $F$
is a $G_\delta$ subset of $X$ and therefore it is also a
continuous image of a Valdivia compact space (see \cite{Ka}). On
the other hand, a well ordered continuous image of a Valdivia
compact space has order type $<\aleph_2$ (see \cite{Ka2}).

It can be shown that $X$ is Valdivia compact if
$\kappa=\aleph_1$. We do not know whether $X$ is Valdivia compact
if $\kappa=\aleph_2$.
\end{uwgi}

\section{A three-space property}

We show that the three-space property for Kadec renormings holds
under the assumption that the quotient space has an LUR renorming.
This solves a problem raised in \cite{LZ} where it is shown that a
Banach space $E$ has a Kadec-Klee renorming provided some subspace
$F$ has a Kadec-Klee renorming and $E/F$ has an LUR renorming.

We begin with an auxiliary lemma on extending Kadec norms.

\begin{lm}\label{wierzba}
Let $E$ be a Banach space and let $F$ be a closed subspace of $E$.
Assume $\Tau$ is a weaker linear topology on $E$ such that $F$ and
the unit ball of $E$ are $\Tau$-closed and $F$ has an equivalent
$\Tau$-Kadec norm. Then there exists an equivalent $\Tau$-lsc norm
$\anorm$ on $E$ which is $\Tau$-Kadec on $F$, i.e. for every $y\in
F$ with $\norm y=1$ and for every $\eps>0$ there exists $V\in
\Tau$ such that $y\in V$ and $\usphere E\cap V\subs \bal(y,\eps)$,
where $\usphere E$ denotes the unit sphere of $E$ with respect to
$\anorm$.
\end{lm}

\begin{pf}
We use ideas from \cite{Ra}. Let $\anorm_0$ be the original norm
of $E$ which, as we may assume, is $\Tau$-lsc and let $B\subs F$
denote the unit closed ball with respect to a given $\Tau$-Kadec
norm. Let $G_n=\cl_\Tau\bal_{\anorm_0}(B,1/n)$. Then each $G_n$ is
a convex, bounded, symmetric neighborhood of the origin in $E$.
Denote by $p_n$ the Minkowski functional of $G_n$ and define
$$\norm x=\sum_{n>0}\al_np_n(x),$$
where $\sett{\al_n}{\ntr}$ is a sequence of positive reals making
the above series convergent. Then $\anorm$ is an equivalent norm
on $E$ which is $\Tau$-lsc, because each $p_n$ is $\Tau$-lsc. We
show that $\anorm$ is $\Tau$-Kadec on $F$.

Fix $y\in F$ with $\norm y=1$ and fix $\eps>0$. By Proposition
\ref{annulus}, find a $\Tau$-neighborhood $W$ of $y$ and $r>1$
such that
$$y\in W\cap rB\subs B(y,\eps/4).$$
We claim that there exist a smaller $\Tau$-neighborhood $U$ of
$y$, $n\in\Nat$ and $\gamma>0$ such that
\begin{equation}
U\cap(r+\gamma)G_n\subs \bal(y,\eps). \tag{1}
\end{equation}
First, find $W_0\in\Tau$ such that $y\in W_0$ and
$W_0+\bal(0,\delta)\subs W$ for some $\delta>0$. Then $W_0\cap
\bal(rB,\delta)\subs \bal(y,\eps/4+\delta)$. Indeed, if $w\in W_0$
and $\norm{w-z}<\delta$ for some $z\in rB$ then $z\in
rB\cap(W_0+\bal(0,\delta))\subs rB\cap W\subs\bal(y,\eps/4)$. Find
$\ntr$ so small that $r/n\loe\delta$ and assume that
$\delta<\eps/4$. Then $W_0\cap\bal(rB,r/n)\subs \bal(y,\eps/2)$.
Next, find $W_1,V\in \Tau$ such that $y\in W_1$, $0\in V=-V$ and
$W_1+V\subs W_0$. Then
$$
W_1\cap\cl_\Tau(\bal(rB,r/n))\subs\bal(y,\eps/2)+V.
$$
Indeed, if $w\in W_1\cap \cl_\Tau(\bal(rB,r/n))$ then there is
$z\in \bal(rB,r/n)$ such that $z-w\in V$, so $z\in W_1+V\subs W_0$
and hence $z\in \bal(y,\eps/2)$. As $V$ can be an arbitrarily
small $\Tau$-neighborhood of $0$, it follows that
$$
W_1\cap \cl_\Tau(\bal(rB,r/n))\subs
\cl_\Tau\bal(y,\eps/2)=\clbal(y,\eps/2).
$$
The last equality follows from the fact that closed balls are
$\Tau$-closed. Note that
$$\cl_\Tau(\bal(rB,r/n))=r\cl_\Tau(\bal(B,1/n))=rG_n.$$ Thus we
have $W_1\cap rG_n\subs\clbal(y,\eps/2)$. Finally, find a
$\Tau$-neighborhood $U$ of $y$ and $\eta>0$ such that
$U+\bal(0,\eta)\subs W_1$ and $\eta<\eps/2$. Let $\gamma>0$ be
such that $\gamma G_n\subs \bal(0,\eta)$. Fix $u\in U\cap
(r+\gamma)G_n$. Then there is $z\in rG_n$ such that $u-z\in \gamma
G_n\subs\bal(0,\eta)$, so $z\in U+\bal(0,\eta)\subs W_1$ and hence
$z\in\clbal(y,\eps/2)$. Thus $u\in
\clbal(y,\eps/2+\eta)\subs\bal(y,\eps)$. This finishes the proof
of (1).

Now, using the fact that each $p_n$ is $\Tau$-continuous on the
$\anorm$-unit sphere, we may assume, shrinking $U$ if necessary,
that $p_n(x)<p_n(y)+\gamma$ whenever $x\in U$ and $\norm x=1$.
Note that $p_n(y)\loe r$, since $r^{-1}y\in G_n$. Thus, if $x\in
U$ and $\norm x=1$ then $p_n(x)<r+\gamma$ which means that $x\in
(r+\gamma)G_n$ and hence $\norm{x-y}<\eps$. This the completes
proof.
\end{pf}

\begin{uwgi}
If, in the above lemma, $\Tau$ is the weak topology then the norm
defined by
$$\norm x = \norm x_0 + \dist(x,F)$$
is Kadec on $F$, where $\anorm_0$ is any equivalent norm such that
$(F,\anorm_0\rest F)$ has the Kadec property. This idea was used
in \cite{LZ}. In general, we do not know whether $\dist(\argum,
F)$ is $\Tau$-lsc.
\end{uwgi}

The following lemma, stated for sequences instead of nets, is due
to Haydon \cite[Proposition~1.2]{Ha} and it is a variation of a
lemma of Troyanski (see \cite[p. 271]{DGZ}) which is an important
tool for obtaining LUR renormings.

\begin{lm}\label{klon}
Let $X$ be topological space, let $S$ be a set and let
$\map{\phi_s,\psi_s}X{[0,+\infty)}$ be lower semi-continuous
functions such that $\sup_{s\in S}(\phi_s(x)+\psi_s(x))<+\infty$
for every $x\in X$. Define
$$
\phi(x)=\sup_{s\in S}\phi_s(x),\qquad \theta_m(x)=
\sup_{s\in S}(\phi_s(x)+2^{-m}\psi_s(x)), \qquad \theta(x)=
\sum_{m\in\nat}2^{-m}\theta_m(x).
$$
Assume further that $\sett{x_\sig}{\Ssig}$ is a net converging to
$x\in X$ and $\theta(x_\sig)\to\theta(x)$. Then there exists a
finer net $\sett{x_\gamma}{\gamma\in\Gam}$ and a net
$\sett{i_\gamma}{\gamma\in\Gam}\subs S$ such that
$$
\lim_{\gamma\in\Gam}\phi_{i_\gamma}(x_{\gamma})=
\lim_{\gamma\in\Gam}\phi_{i_\gamma}(x)=
\lim_{\gamma\in\Gam}\phi(x_{\gamma})=\phi(x)
$$
and
$$
\lim_{\gamma\in\Gam}
(\psi_{i_\gamma}(x_{\gamma})-\psi_{i_\gamma}(x))=0.
$$
\end{lm}

\begin{pf}
By Proposition~\ref{snieg2}, we have
$\lim_{\Ssig}\theta_m(x_\sig)=\theta_m(x)$ for every $m\in\nat$.
Thus, given $m\in\nat$, we can choose $j(m)\in S$ and
$\sig(m)\in\Sig$ such that
$$
\phi_{j(m)}(x)+2^{-m}\psi_{j(m)}(x)>
\sup_{\sig\goe\sig(m)}\theta_m(x_\sig)-2^{-2m}
$$
and
$$
\phi_{j(m)}(x_\sig)>\phi_{j(m)}(x)-2^{-2m}\quad\text{and}\quad
\psi_{j(m)}(x_\sig)>\psi_{j(m)}(x)-2^{-2m}$$ hold for
$\sig\goe\sig(m)$. We may also assume that
$\sig(m_1)\loe\sig(m_2)$ whenever $m_1<m_2$. Define
$$\Gamma=\setof{(\sig,m)\in\Sig\times\nat}{\sig\goe\sig(m)}.$$
Consider $\Gamma$ with the coordinate-wise order and define $\map
h\Gam\Sigma$ by setting $h(\sig,m)=\sig$. Finally, define
$i(\gamma)=j(m)$, where $\gamma=(\sig,m)\in\Gam$. Fix
$\gamma=(\sig,m)\in\Gam$. We have, knowing that $i(\gamma)=j(m)$
and $\sig\goe\sig(m)$,
$$
\phi_{i(\gamma)}(x)+2^{-m}\psi_{i(\gamma)}(x)>
\sup_{\xi\goe\sig(m)}\theta_m(x_\xi)-2^{-2m}\goe
\phi_{i(\gamma)}(x_{h(\gamma)})+2^{-m}\psi_{i(\gamma)}(x_{h(\gamma)})-2^{-2m}.
$$
The last inequality holds because $h(\gamma)=\sig\goe\sig(m)$. It
follows that
$$
|\phi_{i(\gamma)}(x)-\phi_{i(\gamma)}(x_{h(\gamma)})|<2^{-2m+1}\quad
\text{and}\quad
|\psi_{i(\gamma)}(x)-\psi_{i(\gamma)}(x_{h(\gamma)})|<2^{-m+1},
$$
because
$\phi_{i(\gamma)}(x_{h(\gamma)})>\phi_{i(\gamma)}(x)-2^{-2m}$ and
$\psi_{i(\gamma)}(x_{h(\gamma)})>\psi_{i(\gamma)}(x)-2^{-2m}$.
This shows that
\begin{equation}
\lim_{\gamma\in\Gam}|\phi_{i(\gamma)}(x)-\phi_{i(\gamma)}(x_{h(\gamma)})|=0\quad
\text{and}\quad
\lim_{\gamma\in\Gam}|\psi_{i(\gamma)}(x)-\psi_{i(\gamma)}(x_{h(\gamma)})|=0.
\tag{1}\end{equation} We also have
\begin{align*}
\phi_{i(\gamma)}(x)+2^{-m}\sup_{s\in S}\psi_s(x)
&\goe\phi_{i(\gamma)}(x)+2^{-m}\psi_{i(\gamma)}(x)
\goe\sup_{\xi\goe\sig(m)}\theta_m(x_\xi)-2^{-2m}\\
&\goe\limsup_{\eta\in\Gam}\phi(x_{h(\eta)})-2^{-2m}
\goe\liminf_{\eta\in\Gam}\phi(x_{h(\eta)})-2^{-2m}\\
&\goe\phi(x)-2^{-2m} \goe\phi_{i(\gamma)}(x)-2^{-2m}.
\end{align*}
Thus, passing to the limit, we get
\begin{equation}
\liminf_{\gamma\in\Gam}\phi_{i(\gamma)}(x)
\goe\limsup_{\gamma\in\Gamma}\phi(x_{h(\gamma)})
\goe\liminf_{\gamma\in\Gamma}\phi(x_{h(\gamma)}) \goe\phi(x)
\goe\limsup_{\gamma\in\Gam}\phi_{i(\gamma)}(x).
\tag{2}\end{equation} By (1) and (2), the proof is complete.
\end{pf}

\begin{tw}\label{trzy}
Assume $E$ is a Banach space and $\Tau$ is a weaker linear
topology on $X$ such that the unit ball of $E$ is $\Tau$-closed.
Assume further that $F$ is a closed subspace of $E$ which has a
$\Tau$-Kadec renorming and the quotient $E/F$ has a $\Tau'$-lsc
LUR renorming for some Hausdorff locally convex linear topology
$\Tau'$ on $E/F$ such that the quotient map is $(\Tau,\Tau')$
continuous. Then $E$ has a $\Tau$-Kadec renorming.
\end{tw}

Note that since the unit ball of $E/F$ under the LUR renorming is
closed with respect to the weak topology on $E/F$ generated by the
$\Tau'$-continuous linear functionals, we could have equivalently
assumed that $\Tau'$ is the weak topology on $E/F$ generated by a
total subspace of $E/F$.

\begin{pf}
The assumptions imply that $F$ is $\Tau$-closed, being the
pre-image of a singleton under the quotient map. Let $\anorm$ be
an equivalent $\Tau$-lsc norm on $E$ which is $\Tau$-Kadec on $F$
(Lemma~\ref{wierzba}). Denote by $\aabs_q$ the quotient norm on
$E/F$. Let $\aabs$ be an LUR norm on $E/F$ which is $\Tau'$-lsc.
Write $\h x$ for $x+F$, i.e. the image of $x$ under the quotient
map.

Let $\map b{E/F}E$ be a continuous selection for the quotient map
obtained by Bartle-Graves Theorem so that for each $y\in E/F$,
$b(y)\in y$, the range of $b$ on the unit sphere of $E/F$ is
bounded in norm by a positive constant $M$, and $b(ty)=t b(y)$
whenever $t\goe0$ (see \cite[VII Lemma~3.2 and its proof]{DGZ}).
Let $S=\setof{a\in E/F}{\abs a=1}$.

Since the unit ball for $\aabs$ is $\Tau'$-closed, the
$\Tau'$-continuous functionals of unit norm for the dual norm
$\aabs^*$ to $\aabs$ form a norming set for $(E/F,\aabs)$. For
each $a\in S$ choose a $\Tau'$-continuous functional $f_a\in
(E/F)^*$ such that $f_a(a)=1$ and $|f_a|^*\leq2$. Note that if
$\norm{f_a}$ denotes the norm of $f_a$ with respect to $\aabs_q$,
then the values $\norm{f_a}$ are bounded. (We have $|f_a(y)|\leq
|f_a|^*\,|y|\leq 2|y|\leq 2K|y|_q$ for some constant $K$ and hence
$\|f_a\|\leq 2K$.) By enlarging the constant $M$ introduced above,
we may assume that $\|f_a\|\leq M$ for each $a\in S$. Define
$P_ax=f_a(\h x)b(a)$ and let $\psi_a$ be the seminorm given by
$$\psi_a(x)=\norm{x-P_ax}.$$
Note that $\psi_a$ is $\Tau$-lsc, because $P_a$ is a
$\Tau$-continuous functional. Next, define
$$\phi_a(x)=\inf\setof{r>0}{\abs{r^{-1}\h x+a}\loe2}.$$
Observe that $\phi_a$ is the Minkowski functional of the set
$H_a=\setof{x\in E}{\abs{\h x+a}\loe2}$. $H_a$ is a convex set
containing $0$ as an internal point, so $\phi_a$ satisfies the
triangle inequality and is positively homogeneous. Because
$x\mapsto \h x$ is $(\Tau,\Tau')$-continuous and $\aabs$ is
$\Tau'$-lsc, $H_a$ is a $\Tau$-closed set and thus $\phi_a$ is
$\Tau$-lsc. Both families $\setof{\phi_a}{a\in S}$ and
$\setof{\psi_a}{a\in S}$ are pointwise bounded, specifically
$\phi_a(x)\loe\abs{\h x}$ and $\psi_a(x)\loe (M^2+1)\norm x$.
Applying Lemma \ref{klon} we get a $\Tau$-lsc function $\theta$
satisfying the assertion of that lemma and such that $\norm
x_\theta=\theta(x)+\theta(-x)$ defines a $\Tau$-lsc semi-norm on
$E$. Define $\anorm_K$ on $E$ by
$$\norm x_K=\norm x+\abs{\h x}+\norm x_\theta.$$
This is a norm equivalent to $\anorm$. It is $\Tau$-lsc since each
of the three terms defines a $\Tau$-lsc function of $x$. By
Corollary~\ref{snieg1}, the restriction to the unit sphere for
$\anorm_K$ of each of these three functions is $\Tau$-continuous.

We will show that $\anorm_K$ is a $\Tau$-Kadec norm on $E$.

Fix $x\in E$ with $\norm x_K=1$ and fix a net
$\sett{x_\sig}{\Ssig}$ which $\Tau$-converges to $x$ and
$\norm{x_\sig}_K=1$ for every $\Ssig$. We will be done if we find
a finer net converging in norm. We may assume that $x\notin F$, so
that $\h x\not=0$. Since $\|\cdot\|_\theta$ is $\Tau$-continuous
on the sphere, $\lim_{\sig\in\Sigma} \|x_\sig\|_\theta =
\|x\|_\theta$. From the definition of $\|\cdot\|_\theta$ and
Proposition~\ref{snieg2}, we have
$\lim_\Ssig\theta(x_\sig)=\theta(x)$, so by Lemma \ref{klon} we
get a finer net, which we still denote by $\sett{x_\sig}{\Ssig}$
and a net $\sett{a_\sig}{\Ssig}$ such that
\begin{equation}
\lim_\Ssig(\norm{x_\sig-P_{a_\sig}x_\sig}-\norm{x-P_{a_\sig}x})=0
\tag{1}\end{equation}
and
\begin{equation}
\lim_\Ssig\phi_{a_\sig}(x_\sig)=\lim_\Ssig\phi_{a_\sig}(x)=
\lim_\Ssig\sup_{a\in S}\phi_a(x_\sig)=\sup_{a\in S}\phi_a(x).
\tag{2}
\end{equation}
Now observe that $\sup_{a\in S}\phi_a(x)=\abs{\h x}$. Indeed, we
have $\abs{\abs{\h x}^{-1}\h x+a}\loe2$, so $\phi_a(x)\loe\abs{\h
x}$ for every $a\in S$. On the other hand, if $a=\abs{\h x}^{-1}\h
x$, then
$$
\abs{r^{-1}\h x+a}=(r^{-1}+\abs{\h x}^{-1})
\cdot\abs{\h x}=r^{-1}\abs{\h x}+1,
$$
so $\abs{r^{-1}\h x+a}\loe2$ iff $r\goe\abs{\h x}$ which shows
that $\phi_a(x)=\abs{\h x}$.

Let $t=\abs{\h x}^{-1}$.

\begin{claim}\label{sosna} $\lim_{\Ssig}a_\sig=t\h x$. \end{claim}

\begin{pf}
By (2) we have $\lim_\Ssig \phi_{a_\sig}(x)=\abs{\h x}=t^{-1}$.
This means that for every $\eps$ such that $0<\eps<t^{-1}$ there
exists $\sig(\eps)\in\Sig$ such that $\abs{r^{-1}\h x
+a_\sig}\goe2$ whenever $r\loe t^{-1}-\eps$ and
$\sig\goe\sig(\eps)$. Observe that $r^{-1}\h x$ has norm close to
$1$, when $r$ is close to $\abs{\h x}^{-1}$. By LUR, this implies
that $a_\sig$ must be close to $t\h x$. More formally, fix
$\sig\goe\sig(\eps)$ and let $r=t^{-1}-\eps$ and observe that
\begin{align*}
2&\loe\abs{t\h x+a_\sig + (r^{-1}-t)\h x}
\loe \abs{t\h x+a_\sig}+(r^{-1}-t)t^{-1}
= \abs{t \h x+a_\sig}+\Bigl(\frac1{1-\eps t}-1\Bigr).
\end{align*}
It follows that $\liminf_{\Ssig}\abs{t\h x + a_\sig}\goe2$. As
$\abs{t\h x}=1$, the LUR property of $\aabs$ implies
$\lim_{\Ssig}a_\sig=t\h x$.
\end{pf}

By the $(\Tau,\Tau')$-continuity of the quotient map and the
$\Tau$-continuity of $x\mapsto|\h x|$ on the unit sphere, we have
$\Tau'$-$\lim_{\Ssig}\h x_\sig=\h x$ and $\lim_\Ssig\abs{\h
x_\sig}=\abs{\h x}$. As $\aabs$ is a $\Tau'$-lsc LUR norm, it is
$\Tau'$-Kadec and hence by Proposition~\ref{annulus}
$\lim_{\Ssig}\abs{\h x_\sig-\h x}=0$.

\begin{claim}\label{swierk} $\lim_\Ssig P_{a_\sig}x=b(\h x)$.
\end{claim}

\begin{pf} We have
\begin{align*}
t\norm{P_{a_\sig}x-b(\h x)} &=
   \norm{f_{a_\sig}(t\h x)b(a_\sig)-f_{a_\sig}(a_\sig)b(t\h x)}\\
&\loe \norm{f_{a_\sig}(t\h x-a_\sig)
   \cdot b(t\h x)}+\norm{f_{a_\sig}(t\h x)(b(a_\sig)-b(t\h x))}\\
&\loe M\Bigl( \abs{t\h x-a_\sig}_q\cdot\norm{b(t\h x)} +
\norm{tx}\cdot\norm{b(a_\sig)-b(t\h x)} \Bigr).
\end{align*}
By Claim \ref{sosna}, we have $\lim_\Ssig\abs{t\h x-a_\sig}_q=0$
and $\lim_\Ssig\norm{b(a_\sig)-b(t\h x)}=0$ and hence the claim
holds.
\end{pf}

\begin{claim}\label{last}
$\lim_{\Ssig}\norm{P_{a_\sig}x_\sig-P_{a_\sig}x}=0$.
\end{claim}

\begin{pf}
We have
$$
\norm{P_{a_\sig}x_\sig-P_{a_\sig}x}=\norm{f_{a_\sig}(\h x_\sig-\h
x)b(a_\sig)} \loe\norm{f_{a_\sig}}\cdot\abs{\h x_\sig-\h
x}_q\cdot\norm{b(a_\sig)}\loe M^2\abs{\h x_\sig-\h x}_q,
$$
from which the claim follows since $\lim_\Ssig\abs{\h x_\sig-\h
x}=0$ as explained above.
\end{pf}

In order to finish the proof of the theorem, note that (1) and
Claim \ref{swierk} give
$$\lim_\Ssig(\norm{x_\sig-P_{a_\sig}x_\sig}-\norm{x-b(\h x)})=0.$$
Because $\Tau$ is weaker than the norm topology,
Claim~\ref{swierk} and Claim~\ref{last} give $\Tau$-$\lim_\Ssig
P_{a_\sig}x_\sig = b(\h x)$ and hence $\Tau$-$\lim_\Ssig (x_\sig -
P_{a_\sig}x_\sig) = x - b(\h x)$. Thus
$\lim_\Ssig\norm{(x_\sig-P_{a_\sig}x_\sig)-(x-b(\h x))}=0$,
because $\anorm$ is $\Tau$-Kadec on $F$ and $x-b(\h x)\in F$. (If
$x-b(\h x)=0$, use the last displayed equation above instead of
this argument.) Therefore we have
$$
\norm{x_\sig-x}\loe\norm{(x_\sig-P_{a_\sig}x_\sig)-(x-b(\h x))}+
\norm{P_{a_\sig}x_\sig-P_{a_\sig}x}+\norm{P_{a_\sig}x-b(\h x)}.
$$
Since all three of terms on the right tend to $0$, we are done.
\end{pf}

\begin{wn}
Assume $X$ is a locally compact space such that $C_0(X)$ has a
$\Tau_p$-Kadec renorming and $K$ is a compactification of $X$ such
that $C(K\setminus X)$ has a $\Tau_p$-lsc LUR renorming. Then
$C(K)$ has a $\Tau_p$-Kadec renorming. \end{wn}

\begin{pf}
Define $\map T{C(K)}{C(K\setminus X)}$ by setting
$Tf=f\rest(K\setminus X)$. Then $T$ is a bounded, pointwise
continuous linear operator onto $C(K\setminus X)$ and $\ker
T=C_0(X)$. Thus $C(K\setminus X)$ is isomorphic to $C(K)/C_0(X)$.
Apply Theorem \ref{trzy} for $E=C(K)$, $F=C_0(X)$ and $\Tau$,
$\Tau'$ the respective pointwise convergence topologies. \end{pf}

\end{document}